\documentstyle[12pt]{article}
\begin{document}
\newtheorem{lem}{Lemma}[section]
\newtheorem{rem}{Remark}[section]
\newtheorem{question}{Question}[section]
\newtheorem{prop}{Proposition}[section]
\newtheorem{cor}{Corollary}[section]
\newtheorem{thm}{Theorem}[section]
\newtheorem{definition}{Definition}[section]
\newtheorem{exam}{Example}[section]
\newtheorem{conj}{Conjecture}[section]
\newtheorem{openproblem}{Open Problem}[section]
\def\av{{\int \hspace{-2.25ex}-} }
\title
{Some remarks on singular solutions of
 nonlinear elliptic
equations. III:
 viscosity solutions, including parabolic operators
} 
\author{Luis Caffarelli\thanks{Partially
 supported by
        NSF grant DMS-0654267.}\\
Department of Mathematics
\\
The University of Texas
\\
Austin, TX 78712
\\
\\
YanYan Li\thanks{Partially
 supported by
        NSF grant DMS-0701545 and a Rutgers University Research Council 
grant.}\\
Department of Mathematics\\
Rutgers University\\
110 Frelinghuysen Road\\
Piscataway, NJ 08854\\
\\
Louis Nirenberg \\
Courant Institute\\
251 Mercer Street \\
New York, NY 10012\\
}

\date{ }
\maketitle

\input { amssym.def}

\setcounter {section} {-1}

\section{Introduction}

One of the main results in \cite{CLN1},
theorem 1.1, is a strong maximum 
principle for a singular 
supersolution $u$ in a domain $\Omega$ in $\Bbb R^n$
lying above a $C^2$ solution $v$, i.e.  with $u\ge v$.  Recently
we observed that under the conditions in the theorem,   indeed
 under weaker conditions and
also in
 theorems 1.2 and 1.3 there, the function $u$ satisfies in
all of $\Omega$
\begin{equation}
F(x, u, \nabla u, \nabla ^2 u)\le F(x,v,\nabla v, \nabla^2 v)
\quad\mbox{in viscosity sense}.
\label{1}
\end{equation}
(Note that in this paper $F(x, u, \nabla u, \nabla ^2 u)$ amounts
to $F(x, u, \nabla u, - \nabla ^2 u)$ in \cite{CLN1} --- a change of
notation.)

Furthermore, we found that the strong maximum principle holds for functions $u$
which are lower-semi-continuous (LSC) and
satisfy (\ref{1}) in $\Omega$.

Throughout, the nonlinear operator 
$F(x,s,p,M)$ is assumed to be
elliptic for all values of the arguments,
and $C^1$ in $(s,p,M)$, but not uniformly 
elliptic; nor are
$|F_s|$, $
|F_s|$
 uniformly bounded.

In section 1 we prove that the singular functions $u$ 
satisfying the modified condition satisfy (\ref{1}).

We 
would like to point out some new ingredients in the arguments.

In theorem 1.1 in \cite{CLN1} we
considered a function $u$ with possible singularity
at a point, say the origin.  We used a condition that for any $r>0$ small,
\begin{equation}
\inf_{ 0<|x|\le r}
(u+\mbox{any linear function})\ 
\mbox{occurs on}\  \{|x|=r\}.
\label{new2}
\end{equation}
(This condition is related to the notation of
superaffine, as described in \cite{CLN1}.)

In this paper, in section 1 we start by showing
that under a new weaker condition
than (\ref{new2}), a viscosity supersolution
on $0<|x|<r$ of (\ref{1})
becomes a viscosity supersolution in
$|x|< r$.  Namely, we introduce a class of
functions which, to our
knowledge, is new and which
may prove useful in further work:
$\underline{\mbox{lowerconical}}$ functions.  A function $u$ is 
 lowerconical at a point $\bar x\in \Omega$, if
for any $\eta\in C^\infty(\Omega)$, and for any $\epsilon>0$,
$$
\inf_{ x\in \Omega
}
\bigg(  (u+\eta)(x)- (u+\eta)(\bar x) -\epsilon|x-\bar x|\bigg)<0.
$$
This is formulated more precisely in
Definition \ref{definition1}.
Theorem \ref{thm11new}
generalizes Theorem \ref{lem1}
to viscosity supersolutions outside 
a closed manifold,

The notion lowerconical makes sense
also on a Riemannian manifold.
As we pointed out in Remark \ref{remark1.2},
this condition is almost necessary for a viscosity 
supersolution.

We use another ingredient:
a sharpening of the Hopf Lemma.  It is used
in our proof that a viscosity supersolution in a 
punctured region  is also one in the whole region.
The sharp form of the Hopf Lemma, Lemma
\ref{lemA-0new}, refers
 to a linear second order uniformly elliptic
operator $L$ in a bounded 
domain $\Omega$ with $C^2$ boundary,
and to a function $u\ge 0$ in $\Omega$ 
with
$$
u\ge 1\ \ \mbox{in a ball}\ B_\delta\ \mbox{in}\ \Omega.
$$
Lemma
\ref{lemA-0new}
states that there exist $\bar \epsilon$,
$\bar \mu>0$ depending only on 
$n, \delta, \Omega$, and the ellipticity
constants, such that
if
$$
Lu\le \bar \epsilon\ \ \mbox{in}\
\Omega\ \mbox{in viscosity sense},
$$
then 
$$
u(x)\ge \bar \mu \ dist(x, \partial \Omega).
$$

Lemma \ref{lemA-0new} follows from the special case, Lemma \ref{lemA-0},
where $\Omega$ is a ball.
It is a bit surprising that
 we actually need this form of the Hopf Lemma.  In section 1
we give an ``elliptic'' proof of it.  Lemma \ref{lemA-0} is also 
an immediate corollary of a corresponding sharp form
of the Hopf Lemma for parabolic operators,
see Theorem \ref{maximum-parabolic} in
section 5.  Lemma \ref{lemA-0} follows from it by considering $u$ 
independent of $t$.
The analogue of Lemma \ref{lemA-0new},
for parabolic operators, is given in Theorem \ref{thm5.2new}.

Here is an outline of the other sections.  First, section 2 is concerned with 
the maximum principle for LSC viscosity supersolutions $u$ of
 (\ref{1}) in $\Omega$, where $v\in C^2$, in case 
$u\ge v$ on $\partial \Omega$.

\medskip

\noindent{  \underline{Question}.}\
Does the maximum principle hold, i.e., 
$$
u\ge v\quad\mbox{in}\ \Omega,
$$
if, say, $\Omega$ is a small ball?

\medskip 

In general, no, not even for smooth $u$ in case $F$ is not uniformly 
elliptic (see 
Example \ref{example2.1}).  But in Theorem 
\ref{thm3.1}, we prove the maximum principle if
$$
F_u\le 0
$$ --- under a rather mild ellipticity
condition on $F$.

Using a very different kind of argument, in section 1.2, 
we also prove that the maximum principle holds,
without
 assuming $F_u\le 0$, in case $u$ satisfies some linear elliptic inequalities.

In section 3 we prove the strong maximum principle for $u$,
LSC,  satisfying (\ref{1}) in viscosity sense in $\Omega$.  We
also present  an extension of the Hopf Lemma
 for viscosity supersolutions; uniform ellipticity 
is never required.

In section 4 we extend the strong maximum principle 
and the Hopf Lemma to viscosity supersolutions of nonlinear parabolic
operators.  Section 4 is self-contained and may be
read independently of the others.
We thank H. Matano for suggesting that we consider the problem.

\section{Removable singularities for
viscosity solutions}

\subsection{A sufficient condition}
Let $F\in C^0(\Omega\times \Bbb R\times \Bbb R^n\times {\cal S}^{n\times n})$,
where ${\cal S}^{n\times n}$ denotes the set of $n\times n$
real symmetric  matrices and $\Omega$ is a
domain (bounded connected
open set)  in the $n-$dimensional
Euclidean space $\Bbb R^n$.
Throughout the paper we use $B_r(x)$ to denote
a ball  of radius $r$ and centered at $x$,
and use $B_r$ to denote $B_r(0)$.
We use $LSC(\Omega)$ and $USC(\Omega)$ to denote respectively 
the set of lower-semicontinuous and upper-semicontinuous 
functions.

\begin{definition}
Let $\Omega\subset \Bbb R^n$ be an open set, and 
let $u\in LSC(\Omega)$ satisfying
\begin{equation}
\inf_{\Omega}u>-\infty.
\label{BBB3}
\end{equation}
We say that $u$ is lowerconical at a point $\bar x\in \Omega$, if
for any $\eta\in C^\infty(\Omega)$, and for any $\epsilon>0$,
$$
\inf_{ x\in \Omega
} 
\bigg(  (u+\eta)(x)- (u+\eta)(\bar x) -\epsilon|x-\bar x|\bigg)<0.
$$
We say that  $u$ is upperconical
 at $\bar x\in \Omega$, if
 $-u$ is  lowerconical at $\bar x$.

We say that  $u$ is lowerconical on a subset of $E$ of $\Omega$,
if for any $\eta\in C^\infty(\Omega)$, and for any
$\bar x\in E$ and any $\epsilon>0$,
\begin{equation}
\inf_{ x\in \Omega
}
\bigg(  (u+\eta)(x)- (u+\eta)(\bar x) -\epsilon
dist(x, E)\bigg)<0, 
\label{conical}
\end{equation}
where $dist(x,E)$ denotes the distance of $x$ to $E$.
Similarly we say that  $u$ is upperconical
 on $E$ if $-u$ is lowerconical on $E$.
\label{definition1}
\end{definition}

Note that for a smooth submanifold $E$ of dimension
$1\le k\le n-1$, $u(x):=dist(x, E)$ is lowerconical at every point
$\bar x\in E$, but it is not
lowerconical on $E$.

\begin{rem} If $u$ is differentiable at $\bar x$, then $u$ is
both lowerconical and upperconical at  $\bar x$.  In fact, if for some $C^1$
curve $\gamma(t)$ satisfying
$\gamma(0)=\bar x$, $u(\gamma(t))$  is differentiable at 
$0$,  then $u$ is both 
lowerconical  and upperconical at $\bar x$.   On the other hand, 
$u(x)=|x|$, a Lipschitz function, is not lowerconical at $0$.  It is
easy to see that if $\liminf_{x\to \bar x}u(x)>u(\bar x)$, then
$u$ is not lowerconical at $\bar x$.  Also,
$u(x)=\sin(1/|x|)$ for $x\ne 0$, $u(0)=-1$, is both lowerconical
and upperconical at $0$, but is not even continuous.
\end{rem}

\begin{thm}  For $n\ge 1$, let $\Omega$ be a domain
 in $\Bbb R^n$, $\bar x\in \Omega$,
 and  let $F\in C^0(\Omega\times \Bbb R\times
 \Bbb R^n\times {\cal S}^{n\times n})$.
   Assume that   $u\in
LSC(\Omega)$ is lowerconical at 
$\{\bar x\}$ 
and satisfies, for some $f\in USC(\Omega)$,
\begin{equation}
F(x, u, \nabla u, \nabla^2 u)\le f(x),\qquad
\mbox{in}\ \Omega\setminus\{\bar x\}\ \mbox{in
viscosity sense}.
\label{0}
\end{equation}
Then 
\begin{equation}
F(x, u, \nabla u, \nabla^2 u)\le f(x), \qquad
\mbox{in}\ \Omega\ \mbox{in
viscosity sense}.
\label{CCC0}
\end{equation}
\label{lem1}
\end{thm}

\begin{rem} There is a kind of converse.  Namely, 
if $F$ is further assumed to satisfy
\begin{equation}
\limsup_{ a\to \infty} \inf_{  x\in \Omega,
\  |(s,p)|\le \beta  }
F(x, s, p, aI)=\infty,\quad\forall\ \beta>0,
\label{CCC1}
\end{equation} and $f$ is  further assumed to  satisfy
 $\sup_\Omega f<\infty$,
then if
 $u\in LSC(\Omega)$ and satisfies
(\ref{BBB3}) and 
(\ref{CCC0}), necessarily 
 $u$ is lowerconical at every point
of $\Omega$.
On the other hand,  any $u$ satisfies
$\displaystyle{
-e^{-\Delta u}\le 0}$.
This operator does not satisfy (\ref{CCC1}).
  Condition 
(\ref{CCC1}) is clearly  satisfied by uniformly elliptic 
operators.
\label{remark1.2}
\end{rem}

To see the above, suppose that
 $u$  is not lowerconical at some point,
say $0$, in $\Omega$, then for some $\epsilon\in (0,1)$ and some
$\eta\in C^\infty(\Omega)$,
$$
(u+\eta)(x)-(u+\eta)(0)-\epsilon |x|\ge 0,\qquad\mbox{in}\ \Omega.
$$
So for some constant $\delta\in (0, 1)$,
$$
u(x)\ge u(0)-\nabla\eta(0)\cdot x+\frac \epsilon 2 |x|,\qquad
\forall \ 0<|x|<\delta.
$$
For $a>1/(4\delta)$,
$$
u(x)>\varphi_a(x):=u(0)-\nabla\eta(0)\cdot x+a |x|^2,
\qquad \forall\ |x|=\epsilon/(4 a),
$$
$$
\varphi_a(x)\le u(0)+|\nabla \eta(0)|+ \epsilon/(16 a)
\le u(0)+|\nabla \eta(0)|+ 1,\quad
|x|\le \epsilon/(4 a).
$$
Move $\varphi_a$ down, and then up to position $\varphi_a-\bar b$,
$\bar b\ge 0$, so that its graph first touches that of $u$ from below,
at some point $\bar x$, $|\bar x|< \epsilon/(4 a)$.
More precisely,
let
$$
\bar b=\sup \{ b\ |\ 
 u(x)\ge (\varphi_a-b)(x),\ \forall\
|x|\le \epsilon/(4 a)\}.
$$
Clearly, $\bar b\ge 1+u(0)+|\nabla \eta(0)|-\inf_\Omega u$.
On the other hand, since $u(0)=\varphi_a(0)$
and $u(x)>\varphi_a(x)$ for all
$|x|=\epsilon/(4 a)$, we infer that $\bar b\ge 0$, and
for some  $|\bar x|< \epsilon/(4 a)$,
$u(\bar x)=(\varphi_a-\bar b)(\bar x)$ and
$u(x)\ge (\varphi_a-\bar b)(x)$ for all
$|x|\le  \epsilon/(4 a)$.  
By (\ref{CCC0}),
$$
F(\bar x, (\varphi_a-\bar b)(\bar x), \nabla
(\varphi_a-\bar b)(\bar x), \nabla^2 (\varphi_a-\bar b)(\bar x))
\le f(\bar x).
$$
It is easy to see that $|(\varphi_a-\bar b)(\bar x)|$ and
$| \nabla
(\varphi_a-\bar b)(\bar x)|$ are bounded by 
some constant  independent of $a$. Sending $a$ to $\infty$
in the above, we arrive at a contradiction in view of 
(\ref{CCC1}).

\medskip

The following example shows that the
assumption on $u$ in Theorem \ref{lem1} is essentially optimal.

\begin{exam}  Let $u(x)=|x|$.
Then
$$
F(x,u, \nabla u, \nabla^2 u):= -e^{-\Delta u}+1-|\nabla u|^2\le 0,
\qquad
\mbox{in}\ B_1\setminus\{0\}.
$$
But
the inequality does not hold in $B_1$ in viscosity sense, as easily
seen by taking $\varphi(x)=|x|^2$ as a test function.
\end{exam}

\noindent{\bf Proof of Theorem \ref{lem1}.}\  We may assume that 
$\bar x=0$. 
Let $\varphi\in C^2(\Omega)$ satisfy
$(u-\varphi)(0)=0$, $u-\varphi\ge 0$ in 
$\Omega$.
For any $0<\delta< dist(0, \partial \Omega)/9$, we consider
\begin{equation}
\varphi_\delta(x):= \varphi(x)-\frac \delta 2 |x|^2.
\label{FF1}
\end{equation}
Clearly,
\begin{equation}
u(x)>\varphi_\delta(x),\qquad  x\in \Omega\setminus\{0\},
\label{2}
\end{equation}
\begin{equation}
u(x)\ge \varphi_\delta(x) +\frac 12 \delta^3,\qquad   x\in \Omega,\  
|x|\ge \delta.
\label{3}
\end{equation}

Since $u\in LSC(\Omega)$ is lowerconical, we have, for large $i$, 
$$
\liminf_{ x\to 0}
\bigg(  (u-\varphi_\delta)(x)- (u-\varphi_\delta)(0) -
\frac 1i |x|\bigg)
\ge 0,
$$
and 
there exists $\{x_i\}\subset \Omega\setminus \{
0\}$ such that

\begin{eqnarray}&&(u-\varphi_\delta)(x_i)- \frac 1i |x_i|
\nonumber\\
&=&(u-\varphi_\delta)(x_i)- (u-\varphi_\delta)(0) - \frac 1i |x_i|
\nonumber\\
&=&\inf_{ x\in \Omega}
\bigg(  (u-\varphi_\delta)(x)- 
(u-\varphi_\delta)(0) - \frac 1i|x|\bigg)<0.
\label{4}
\end{eqnarray}

Claim that 
\begin{equation}
\lim_{i\to \infty} x_i=0.
\label{5}
\end{equation}

Indeed, let $x_i\to \hat x$ along a subsequence, still
denoted as $\{x_i\}$. 
Then, after sending $i$ to infinity in (\ref{4}),
we have
$$
(u-\varphi_\delta)(\hat x)\le 0,
$$
which implies 
$\hat x=0$ in view of (\ref{2}) and (\ref{3}).
We have proved (\ref{5}).

\vskip 5pt
\hfill $\Box$
\vskip 5pt

Let
$$
\varphi_\delta^{(i)}(x):=
\varphi_\delta(x)+ \frac 1{\sqrt{i}}
\frac{ x_i }{  |x_i| }\cdot x.
$$
We have, in view of (\ref{4}) and (\ref{5}), that
$$
(u-\varphi_\delta^{(i)})(x_i)
= (u-\varphi_\delta)(x_i)
- \frac 1{\sqrt{i}}  |x_i| 
< (\frac 1i -\frac 1{\sqrt{i}}) |x_i| <0.
$$

Since 
$$
(u-\varphi_\delta^{(i)})(0)=0,
$$
and, in view of (\ref{2}) and (\ref{3}),
$$
(u-\varphi_\delta^{(i)})(x) \ge \frac 12 \delta^3 +O(\frac 1{\sqrt{i}})>0,\qquad
\forall\ |x|\ge \delta, \ \ \mbox{for large}\ i,
$$
there exists $\widetilde x_i$, $0<|\widetilde x_i|<\delta$,
such that
$$
(u-\varphi_\delta^{(i)})(\widetilde x_i)=
\min_{ 0<|x|<\delta } (u-\varphi_\delta^{(i)})(x)<0.
$$
Namely,
$$
\psi_\delta^{(i)}(x):= \varphi_\delta^{(i)}(x)+ (u-\varphi_\delta^{(i)})(\widetilde x_i)
$$
satisfies
$\psi_\delta^{(i)}(\widetilde x_i)=u(\widetilde x_i)$ and
$u\ge \psi_\delta^{(i)}$ near $\widetilde x_i$.

Similar to (\ref{5}), we have
$$
\lim_{i\to \infty} \widetilde x_i=0.
$$

Thus, by 
(\ref{0}),
$$
F(\widetilde x_i, \psi_\delta^{(i)}(\widetilde x_i),
\nabla \psi_\delta^{(i)}(\widetilde x_i), \nabla^2 \psi_\delta^{(i)}(\widetilde x_i))
\le f(\widetilde x_i).
$$
Sending $i$ to $\infty$ in the above leads to
$$
F(0, \varphi_\delta(0), \nabla \varphi_\delta(0), \nabla^2 
\varphi_\delta(0))\le f(0).
$$
Sending $\delta$ to $0$ in the above leads to
$$
F(0, \varphi(0), \nabla \varphi(0), \nabla^2 
\varphi(0))\le f(0).
$$
Theorem \ref{lem1} is proved.

\vskip 5pt
\hfill $\Box$
\vskip 5pt

\subsection{Another sufficient condition for
removable singularity}

Let
$(a_{ij}(x))$,  $b_i(x)$ and $c(x)$
be  $L^\infty(\Omega)$ functions satisfying,
for some positive constants $\lambda$ and $\Lambda$,
\begin{equation}
|a_{ij}(x)|+|b_i(x)| +|c(x)|\le \Lambda,
\ \
a_{ij}(x)\xi_i\xi_j\ge \lambda|\xi|^2,\qquad\forall\
x\in \Omega, \xi\in \Bbb R^n.
\label{aij}
\end{equation}

In the rest of this section we assume that $F$
 is a degenerate elliptic operator:
\begin{equation}
F(x, s, p, M+N)\ge F(x, s, p, M),\qquad  \forall\ (x, s, p, M)\in
\ \ N\in {\cal S}^{n\times n}_+,
\label{elliptic}
\end{equation}
where ${\cal S}^{n\times n}_+\subset {\cal S}^{n\times n}$ denotes the set of
positive definite matrices.

\begin{thm}
For $n\ge 1$, let $F\in C^0(\Omega\times \Bbb R\times
 \Bbb R^n\times {\cal S}^{n\times n})$ satisfy
(\ref{elliptic}),   $(a_{ij}(x))$,  $b_i(x)$ and $c(x)$
be as above, and let $f\in USC(\Omega)$.
Assume that $u\in
LSC(\Omega)$ satisfies, for some constant $C$,
\begin{equation}
a_{ij}(x)\partial_{ij}u+b_i(x)\partial_iu+c(x)u\le C,
\quad\mbox{in}\ \Omega, \
\mbox{in viscosity sense,}
\label{B1-newnew}
\end{equation}
and, for some subset $E$ of $\Omega$ of zero Lebesgue measure,
\begin{equation}
F(x, u, \nabla u, \nabla^2 u)\le f(x) \qquad
\mbox{in}\ \Omega\setminus  E\  \mbox{in the viscosity sense}.
\label{abab5}
\end{equation}
 Then
$$
F(x, u, \nabla u, \nabla^2 u)\le f(x) \qquad
\mbox{in}\ \Omega\  \mbox{in the viscosity sense}.
$$
\label{thmnew1}
\end{thm}

\noindent{\bf Proof of Theorem \ref{thmnew1}.}\
Let $\varphi\in C^2(\Omega)$ satisfy
$$
\varphi \le u,\quad\mbox{in}\ \Omega,\qquad \mbox{and}\ \varphi(\bar
x)=u(\bar x),\ \mbox{for some}\
 \bar x\in \Omega.
$$
We have only to prove that
\begin{equation}
F(\bar x, \varphi(\bar x),  \nabla  \varphi(\bar x),  \nabla^2
\varphi(\bar x))\le f(\bar x).
\label{aa1}
\end{equation}
We need only consider $\bar x\in E$ and
may assume,
without  loss of generality, that
$\bar x=0\in E$.
 For any $0<\delta<dist(0, \partial \Omega)/9$,
let $\varphi_\delta$ be defined in (\ref{FF1}).
Then (\ref{2}) and (\ref{3}) hold.  

For $\epsilon\in (0, \delta^3/4)$, 
let,
$$
w_\epsilon(x)\equiv w^{(\delta)}_\epsilon(x):=
\left\{
\begin{array}{ll}
\min\{(u-\varphi_\delta)(x)-\epsilon, 0\}& x\in B_\delta,\\
0& x\in B_{2\delta}\setminus B_\delta,
\end{array}
\right.
$$
and
$$
\Gamma_{w_\epsilon} (x):=
\sup\{a+b\cdot x\ |\ a\in \Bbb R, b\in \Bbb R^n,
a+b\cdot z\le w_\epsilon(z)\ \forall\ z\in
B_{2\delta}\}
$$
be the convex envelope of $w_\epsilon$ on $B_{2\delta}
\equiv B_{2\delta}(0)$.

Since $w_\epsilon=0$ outside $B_\delta$,
 and
$\min_{ B_{2\delta} }\Gamma_{w_\epsilon}\le w_\epsilon(0)=
-\epsilon<0$,
the contact set of $w_\epsilon$ and $\Gamma_{w_\epsilon}$ satisfies
\begin{equation}
\{ x\in B_{2\delta}\
|\  w_\epsilon(x)=\Gamma_{w_\epsilon}(x) \}
\subset \{ x\in B_\delta |\  w_\epsilon(x)=
(u-\varphi_\delta)(x)-\epsilon<0 \}.
\label{G2}
\end{equation}

We will need
\begin{prop}  There exists some positive constants $K$
 such that
for any point $\bar x\in \{ x\in B_{2\delta}\
|\  w_\epsilon(x)=\Gamma_{w_\epsilon}(x) \}$,
there exists $\bar p\in \Bbb R^n$ so that
$$
\Gamma_{w_\epsilon}(x) \le
\Gamma_{w_\epsilon}(\bar x)+\bar p\cdot (x-\bar x)+
K|x-\bar x|^2, \qquad \forall\ x\in B_{2\delta}.
$$
\label{lem-35}
\end{prop}

The proof of this proposition will be postponed to the end of the
proof of the theorem.

Once Proposition \ref{lem-35} is proved, we can apply, as in section 3 of
\cite{CLN1}, 
lemma 3.5 of
\cite{CC} to obtain
 that $\Gamma_{w_\epsilon}\in C^{1,1}_{loc}(B_{2\delta})$, and
then use
 the Alexandrov-Bakelman-Pucci estimate to
obtain
$$
\epsilon^n=|\inf_ {  B_{2\delta} } w_\epsilon|^n
\le \int_{ \{w_\epsilon=\Gamma_{w_\epsilon} \} } \det(\nabla^2
\Gamma_{w_\epsilon}).
$$
This implies that
$$
\mbox {The Lebesgue measure of}\  \{w_\epsilon=\Gamma_{w_\epsilon} \} >0.
$$

Since $\Gamma_{w_\epsilon}$ is convex, 
it is, by the 
 Alexandrov theorem,
  second order 
differentiable except on a set of
zero Lebesgue measure.
We also know that $E$ has zero Lebesgue measure.

Thus we can pick a point $x$ ($=x_\epsilon$)
in  $ \{w_\epsilon=\Gamma_{w_\epsilon} \}
\cap (B_{2\delta} \setminus E)$
where $\Gamma_{w_\epsilon}$ is 
second order
differentiable.

We know from (\ref{G2}) and the definition of 
$\Gamma_{w_\epsilon}$ that
\begin{equation}
u(x)= \psi_\epsilon(x)
:= \varphi_\delta(x)+\epsilon+\Gamma_{ w_\epsilon}(x),
\label{bb0}
\end{equation}
and
$$
u\ge \psi_\epsilon=   \varphi_\delta+\epsilon+\Gamma_{ w_\epsilon},
\qquad \mbox{near}\ x.
$$
Since $ \Gamma_{w_\epsilon}$ is second order
differentiable at $x$, $\nabla \Gamma_{w_\epsilon}(x)$ is well
defined, and,
for any $\mu>0$, and for $z$ near $x$,
\begin{eqnarray*}
u(z)&\ge &
(\Gamma_{w_\epsilon}+\varphi_\delta)(x)+\epsilon+
\nabla(\Gamma_{w_\epsilon}+\varphi_\delta)(x)\cdot (z-x)\\
&&  +\frac 12 (z-x)^t \nabla^2(\Gamma_{w_\epsilon}+\varphi_\delta)(x)
\cdot (z-x)+\circ(|z-x|^2)\\
&\ge &
(\Gamma_{w_\epsilon}+\varphi_\delta)(x)+\epsilon+
\nabla(\Gamma_{w_\epsilon}+\varphi_\delta)(x)\cdot (z-x)\\
&&  +\frac 12 (z-x)^t \nabla^2(\Gamma_{w_\epsilon}+\varphi_\delta)(x)
\cdot (z-x)-\frac {\mu}2 |z-x|^2.
\end{eqnarray*}

By (\ref{abab5}) and the above,
$$
F(x, (\Gamma_{w_\epsilon}+\varphi_\delta)(x)+\epsilon, \nabla 
(\Gamma_{w_\epsilon}+\varphi_\delta)(x), \nabla^2 
(\Gamma_{w_\epsilon}+\varphi_\delta)(x)-\mu I)\le f(x).
$$ 

Sending $\mu$ to $0$ leads to
\begin{equation}
F(x, (\Gamma_{w_\epsilon}+\varphi)(x)+\epsilon, \nabla
(\Gamma_{w_\epsilon}+\varphi)(x), \nabla^2
(\Gamma_{w_\epsilon}+\varphi)(x))\le f(x).
\label{bb2}
\end{equation}
Clearly,
\begin{equation}
|\Gamma_{w_\epsilon}(x)|\le  \epsilon,
 \qquad |\nabla \Gamma_{w_\epsilon}(x)|\le \frac
 \epsilon \delta.
\label{bb1}
\end{equation}
By the convexity of $\Gamma_{w_\epsilon}$, $\nabla^2 \Gamma_{w_\epsilon}(x)
\ge 0$.
We see from (\ref{bb0}) and (\ref{bb1}) that (recall that
$x=x_\epsilon$)
$ (u-\varphi_\delta)(x_\epsilon)\to 0$ as $\epsilon\to 0$.
This 
which implies, in view of (\ref{2}),
(\ref{elliptic}), (\ref{bb2}) and the convexity of $\Gamma_{w_\epsilon}$,
 that
$
 x_\epsilon\to 0
$
and
\begin{eqnarray}
f(0)&\ge& 
\limsup_{\epsilon\to 0}
F(x, (\Gamma_{w_\epsilon}+\varphi_\delta)(x)+\epsilon, \nabla
(\Gamma_{w_\epsilon}+\varphi_\delta)(x), \nabla^2
(\Gamma_{w_\epsilon}+\varphi_\delta)(x))\nonumber\\
&\ge & \limsup_{\epsilon\to 0}
F(x, (\Gamma_{w_\epsilon}+\varphi_\delta)(x)+\epsilon, \nabla
(\Gamma_{w_\epsilon}+\varphi_\delta)(x), \nabla^2
\varphi_\delta)(x))\nonumber\\
&=&F(0, \varphi_\delta(0),
\nabla
 \varphi_\delta(0),
 \nabla^2
 \varphi_\delta(0))
\nonumber\\
&=& 
F(0, \varphi(0),
\nabla
 \varphi(0),
 \nabla^2
 \varphi(0)-\delta I)\nonumber.
\end{eqnarray}
Sending $\delta$ to $0$ in the above leads to 
$$
F(0, \varphi(0),
\nabla
 \varphi(0),
 \nabla^2
 \varphi(0)\le f(0).
$$
Theorem \ref{thmnew1} is established --- provided Proposition \ref{lem-35}
holds.

\bigskip

Now we prove Proposition \ref{lem-35}.
  Under $\Delta u\le C$ instead of (\ref{aij}), 
the above proposition was proved in \cite{CLN1},
see lemma 3.1 there.
The new ingredient is the following

\subsection{A strengthening of the Hopf Lemma}

\begin{lem}
Let $\Omega$ be a domain of $\Bbb R^n$, with
$C^2$ boundary,
 and let  $(a_{ij}(x))$,  $b_i(x)$ and $c(x)$
be  $L^\infty(\Omega)$  functions satisfying
(\ref{aij})
for some positive constants $\lambda$ and $\Lambda$.
 Let $B\subset \Omega$ be a ball of radius
$\delta$.
Then  there exist some positive constants $\bar \epsilon, \bar \mu>0$
which depend only on $n, \lambda, \Lambda, \delta,
\Omega$
such that
if
  $u\in
LSC(\Omega)$   satisfies
$$
a_{ij}(x)\partial_{ij}u+b_i(x)\partial_iu+c(x)u\le  \bar \epsilon,
\quad \mbox{in}\ \Omega, \
\mbox{in viscosity sense},
$$
$$
u\ge 0,\ \ \mbox{in}\ \Omega, \quad \mbox{and}\ \ 
u\ge 1,\qquad \mbox{on}\ B.
$$
Then
$$
u(x)\ge \bar \mu
\ dist(x, \partial \Omega),\quad \mbox{in}\ \Omega.
$$
\label{lemA-0new}
\end{lem}

We first prove Lemma  \ref{lemA-0new} for $\Omega=B_1$, which is stated as

\begin{lem}
Let
  $(a_{ij}(x))$,  $b_i(x)$ and $c(x)$
be  $L^\infty(\Omega)$  functions satisfying
(\ref{aij}) with $\Omega=B_1$
for some positive constants $\lambda$ and $\Lambda$.
Then, for any $0<\delta<1$, there exist
some positive constants $\bar \epsilon, \bar \mu>0$
which depend only on $n, \lambda, \Lambda, \delta$,
such that
if 
  $u\in
LSC(\Omega)$   satisfies
\begin{equation}
Lu:= a_{ij}(x)\partial_{ij}u+b_i(x)\partial_iu+c(x)u\le  \bar \epsilon,
\quad \mbox{in}\ B_1, \
\mbox{in viscosity sense},
\label{A1}
\end{equation}
\begin{equation}
u\ge 1,\qquad \mbox{on}\ B_{\delta}\subset B_1,
\label{A2}
\end{equation}
and
$$
u\ge 0,
\quad \mbox{in}\ B_1.
$$
Then
$$
u\ge \bar \mu (1-|x|),\qquad \mbox{on}\ B_1.
$$
\label{lemA-0}
\end{lem}
 
Lemma  \ref{lemA-0new} then follows by  repeated application
of this for scaled balls.

\noindent{\bf Proof of Lemma \ref{lemA-0}.}\ 
For a large positive constant $\alpha$ to be chosen later, 
consider the function
$$
v(x, t):= \frac {u(x)}{  \cos(\alpha t) },\quad
\mbox{in}\ B_1\times (-\beta, \beta),\ \ 
\beta:= \frac \pi {10 \alpha}.
$$
In particular, consider $v$ in the ellipsoid
$$
E_1:=\{(x, t)\ |\ |x|^2+\beta t^2<1\}.
$$
By (\ref{A2}),
$$
v\ge 1,\qquad\mbox{in}\ 
E_\delta:= \{(x, t)\ |\ |x|^2+\beta t^2<
\delta^2\}.
$$
By (\ref{A1}), 
$$
(L+\partial_t^2)u\le \bar \epsilon, \quad \mbox{in}\ E_1, \
\mbox{in viscosity sense}.
$$
A computation gives
$$
(L+\partial_t^2)u=
\cos(\alpha t)( Lv +v_{tt})
-2\alpha \sin (\alpha t) v_t -\alpha^2 (\cos(\alpha t)) v.
$$
Hence
$$
\widetilde L v: =
a_{ij}v_{ij}+b_iv_i+v_{tt}
-2\alpha\tan (\alpha t)v_t- (\alpha^2-c)v
\le \frac {\bar\epsilon}{\cos(\alpha t)}\le 2\bar \epsilon.
$$
We now fix the value of $\alpha$ to be $\sqrt{\Lambda}$.
Then,
\begin{equation}
 (\alpha^2-c) \ge 0,\qquad \mbox{in}\ 
E_1.
\label{A7}
\end{equation}
In $E_1\setminus E_\delta$, consider the comparison function
$$
h(x,t):=\frac 
{E-e^{-k}}
D, \quad E:= e^{ -k(|x|^2+\beta t^2) },\quad 
\qquad D:= e^{ -k \delta^2}-e^{-k}.
$$
Then
$$
h_i=-2k x_i \frac ED, \quad  h_t= -2k\beta t  \frac ED,
$$
$$
h_{ij}=(4k^2 x_ix_j-2k\delta_{ij}) \frac ED,
\quad h_{tt}=(4k^2\beta^2 t^2-2k\beta) \frac ED.
$$
Hence, for any constant $a\ge 0$,
\begin{eqnarray*}
\widetilde L(h-a)&=&
\frac ED\bigg\{
a_{ij}(4k^2 x_ix_j-2k\delta_{ij})
-b_i(2kx_i)+
(4k^2\beta^2 t^2 -2k \beta)\\
&&+2\alpha \tan(\alpha t)(2k\beta t)
-(\alpha^2 -c)\bigg\}
+(\alpha^2-c)\frac {e^{-k}}D
+a(\alpha^2 -c).
\end{eqnarray*}

Now move $h$ down, and then up to position $h-a$, $a\ge 0$, so that
its graph first touches that of $v$ from below, at some point
$(\bar x, \bar t)$.  We claim that
$a=0$, so that 
$
u\ge h$
and we would have the desired conclusion.
To see this, suppose $a>0$, then $(\bar x, \bar t)\in
E_1\setminus \overline E_\delta$.
Thus, in view of (\ref{A7}), we have at 
 $(\bar x, \bar t)$ that
\begin{eqnarray*}
2\bar \epsilon & \ge &  \widetilde Lv \ge  \widetilde L(h-a)
\ge \frac ED 
\bigg\{
a_{ij}(4k^2 x_ix_j-2k\delta_{ij})
-b_i(2kx_i)\\
&& +
(4k^2\beta^2 t^2 -2k \beta)+2\alpha \tan(\alpha t)(2k\beta t)
-(\alpha^2 -c)\bigg\}.
\end{eqnarray*}
It follows that for some constants $\bar k, c_0>0$, depending only on $ 
n, \lambda, \Lambda$ and $\delta$,
$$
2\bar \epsilon \ge   
 c_0 \bar k^2\frac ED.
$$
Then at $(\bar x, \bar t)$ we have
$$
2\bar \epsilon 
(e^{ -\bar k \delta }
-e^{-\bar k}) e^{  \bar k(|x|^2 +\beta t^2)}\ge
c_0\bar k^2.
$$
This 
 shows that
if $\bar \epsilon$ is small then this is impossible.

\vskip 5pt
\hfill $\Box$
\vskip 5pt

The following example shows that the smallness of $\bar \epsilon$
in Lemma \ref{lemA-0}  indeed
depends on $\delta$.

 In $B_1\subset \Bbb R^2$, consider, for $0<\delta<1/4$, 
the function
$$
w(x)
=
\left\{
\begin{array}{ll}
-\log(|x|+\delta)(1-\delta-|x|)^2,
&\qquad \mbox{for}\ 0\le |x|\le 1-\delta, \\
0,&\qquad \mbox{for}\ 1-\delta\le |x|\le 1.
\end{array}
\right.
$$
It is easy to check that
 $w\in C^2(\overline B_1)$, and
$$
\Delta w\le C,\quad\mbox{for some constant
independent of}\ \delta.
$$
Now 
$$
w(|x|)\ge -\log(2\delta) (1-2\delta)^2,\qquad |x|\le \delta.
$$
Hence
$$
u:= w/[  -\log(2\delta) (1-2\delta)^2 ]
$$
satisfies
$$
u\ge 0\ \ \mbox{in}\
B_1,\qquad
u\ge 1\ \ \mbox{in}\
B_\delta,
$$
and, for some positive constant $C'$ independent of $\delta$,
$$
Lu\le  \frac{C'}{ |\log \delta|}.
$$
But 
$u(x)=0$ for $1-\delta\le |x|\le 1$.  
Thus in the lemma, it is necessary
that  $\bar \epsilon < C''/|\log\delta|$ 
for some $C''$ smaller than $C'$.

We will also use

\begin{lem} 
  Let 
 $(a_{ij}(x))$,  $b_i(x)$ and $c(x)$
satisfy (\ref{aij}) for some  positive constants $\lambda$ and $\Lambda$,
$u\in LSC(\Omega)$ satisfy (\ref{B1-newnew}),
 and
let   $\omega$ be a
 non-negative  non-decreasing continuous function  on
$(0, 2d), \ d:=diam(\Omega))$.
Assume that   $u$ satisfies,
for some $\bar x, \bar y\in \Omega$,
$\bar x \ne \bar y$,  and  $p, q$  in $\Bbb R^n$,
that
$$
u(y)\ge u(\bar x)+p\cdot (y-\bar x)-|y-\bar x|\omega(|y-\bar x|),\qquad
y\in \Omega,
$$
and
$$
u(z)\ge u(\bar x)+p\cdot (\bar y-\bar x)-
|\bar y-\bar x|\omega(|\bar y-\bar x|)+ q\cdot (z-\bar y)
-|z-\bar y|\omega(|z-\bar y|),\ \
\forall\ z\in \Omega.
$$
Then,  for some positive constants $C_1$ and $C_2$ depending only
on $n$, $\lambda$, $\Lambda$, $C$, 
$$
|p-q|\le C_1
\omega(2|\bar x-\bar y|)
+C_2 C|\bar x-\bar y|,
$$
\label{lem-F3}
\end{lem}
where $C$ is the constant in (\ref{B1-newnew}).

\noindent{\bf Proof of Lemma \ref{lem-F3}.}\
Working with $\tilde u(z)= u(z+\bar x)-[u(\bar x)+p\cdot z]$
instead of $u(z)$, we may assume without loss of generality
that $\bar x=0$, $u(0)=0$, $p=0$, $\bar y\ne 0$, $q\ne 0$:
\begin{equation}
u(z)\ge -|z|\omega(|z|), \qquad z\in \Omega,
\label{F4}
\end{equation}
\begin{equation}
u(z)\ge - |\bar y|
\omega(|\bar y|)+ q\cdot (z-\bar y)
-|z-\bar y|\omega(|z-\bar y|),\qquad
\forall\ z\in \Omega.
\label{F5}
\end{equation}
By (\ref{F4}),
$$
u(z)\ge - 2|\bar y| \omega(2|\bar y|),
\qquad \forall\ 0<|z|\le 2|\bar y|.
$$
For  $|z-\bar y|\le \frac 12 |\bar y|$,
 we deduce from
 (\ref{F5}) that
$$
u(z)\ge
 -
|\bar y|\omega(|\bar y|)
+ q\cdot (z-\bar y)
-\frac 12 |\bar y|
\omega(\frac 12|\bar y|)
\ge
-
2|\bar y|\omega(2|\bar y|)
+ q\cdot (z-\bar y).
$$
It follows that
\begin{equation}
u(z)\ge
 -
2|\bar y|\omega(2|\bar y|) +\frac 14 |q||\bar y|,
\qquad \forall\ z\in
(B_{ 2 |\bar y|}(0) \setminus   B_{ |\bar y|/2}(\bar y))\cap U,
\label{F9}
\end{equation}
where
$$
U:= \{z\in \Bbb R^n\ |\
q\cdot (z-\bar y)\ge \frac 12 |q||z-\bar y|\}.
$$

We may assume that
$$
|q|\ge 32 \omega(2|\bar y|),
$$
since otherwise there is nothing to prove.
So we deduce from (\ref{F9}) that
$$
u(z)\ge
\frac 18  |q||\bar y|\ge 4|\bar y|\omega(2|\bar y|),
\qquad \forall\ z\in
(B_{ 2 |\bar y|}(0) \setminus   B_{ |\bar y|/2}(\bar y))\cap U,
$$

Let
$$
\widetilde u(x):= u(|\bar y|x)/(|q||\bar y|),
\qquad
x\in B_2.
$$
Then we have
$$
\widetilde u(0)=0,
$$
$$
\widetilde u \ge -2\omega(2|\bar y|)/|q|,
\qquad
\mbox{in}\ B_2,
$$
$$
\widetilde u \ge \frac 
18, \qquad (B_2  \setminus  B_{1/2}(e))\cap \widetilde U,
$$
and
$$
\widetilde a_{ij}\partial_{ij}\widetilde u
+|\bar y| \widetilde b_i\partial _i \widetilde u
+|\bar y|^2 \widetilde c\widetilde u
\le C|\bar y|/|q|,\qquad
\mbox{in}\ B_2,
$$
where $e=\bar y/|\bar y|$,
$$
\widetilde U:= \{x\in \Bbb R^n\ |\
q\cdot (x-e)\ge \frac 12 |q||x-e|\},
$$
and
$$
\widetilde a_{ij}(x)=a_{ij}(|\bar y|x),\quad
\widetilde b_i(x)=b_i(|\bar y|x),\quad
\widetilde c(x)=c(|\bar y|x).
$$

Applying Lemma 
\ref{lemA-0} to
$\displaystyle{
u(x)= 8[\widetilde u(2x)+2\omega(2|\bar y|)/|q|]
}$,
 we have, for some $\bar \epsilon>0$,
depending only on $n, \lambda, \Lambda$,
$$
\frac {\omega(2|\bar y|)}{ |q| }>\bar \epsilon,
\quad \mbox{or}\quad   
\frac{C |\bar y| }{ |q| }>\bar\epsilon.
$$
The desired estimate follows.
Lemma \ref{lem-F3} is established.

\vskip 5pt
\hfill $\Box$
\vskip 5pt

\noindent{\bf Proof of Proposition \ref{lem-35}.}\
Given Lemma \ref{lem-F3}, the proof is the same as that of
lemma 3.1 in 
\cite{CLN1} --- using this lemma instead of 
lemma A there.

\vskip 5pt
\hfill $\Box$
\vskip 5pt

The proof of Theorem  \ref{thmnew1} is completed.

\subsection{A supersolution outside a closed submanifold}

\begin{thm}  For $n\ge 1$,  let $F\in
C^0(\Omega\times \Bbb R\times
 \Bbb R^n\times {\cal S}^{n\times n})$ satisfy
(\ref{elliptic}), and   let
$\Omega\subset \Bbb R^n$ be a domain,
 $E\subset \Omega$ be a smooth closed
 submanifold of dimension $k$,  $0\le k\le n-1$,
 and let $f\in USC(\Omega)$.
Assume that     $u\in LSC(\Omega)$
is lowerconical in $E$, and
satisfies
\begin{equation}
F(x, u, \nabla u, \nabla^2 u)\le f(x) \qquad
\mbox{in}\ \Omega\setminus E,\ \ \mbox{in the viscosity sense}.
 \label{abcde} \end{equation} Then
$$
F(x, u, \nabla u, \nabla^2 u)\le f(x) \qquad
\mbox{in}\ \Omega, \ \mbox{in the viscosity sense}.
$$
\label{thm111new}
\end{thm}

\begin{thm}  For $n\ge 2$, let $F\in
C^0(\Omega\times \Bbb R\times
 \Bbb R^n\times {\cal S}^{n\times n})$ satisfy
(\ref{elliptic}), and   let 
$\Omega\subset \Bbb R^n$ be a domain,
 $E\subset \Omega$ be a smooth closed
 submanifold of dimension $k$,  $0\le k\le n-2$, and let $f\in 
USC(\Omega)$.
  Assume that   $u\in
LSC(\Omega\setminus E)$  satisfies
\begin{equation}
\inf_{ \Omega\setminus E} u>-\infty,
\label{ABC0}
\end{equation}
and
\begin{equation}
\lambda_1(\nabla^2 u)+\cdots+\lambda_{k+2}(\nabla^2 u)\le 0,
\qquad \mbox{in}\
\Omega\setminus E, \ \mbox{in the viscosity sense},
\label{ABC1}
\end{equation}
Then, after extending $u$ to $E$ by letting
$$
u(x):= \liminf_{y\in \Omega\setminus E, \  y\to x}
u(y),
$$
$u$ is  in
$LSC(\Omega)$, and is lowerconical in $E$.
\label{thm11new}
\end{thm}

A corollary of the above two theorems is

\begin{cor}  For $n\ge 2$, let $F\in
C^0(\Omega\times \Bbb R\times
 \Bbb R^n\times {\cal S}^{n\times n})$ satisfy
(\ref{elliptic}), and   let
$\Omega\subset \Bbb R^n$ be a domain,
 $E\subset \Omega$ be a smooth closed
 submanifold of dimension $k$,  $0\le k\le n-2$, and let $f\in
USC(\Omega)$.
  Assume that   $u\in
LSC(\Omega\setminus E)$  satisfies (\ref{ABC0}), 
(\ref{ABC1}) and
$$
F(x, u, \nabla u, \nabla^2 u)\le f(x) \qquad
\mbox{in}\ \Omega\setminus E, \ \mbox{in the viscosity sense}.
$$
Then
$$
F(x, u, \nabla u, \nabla^2 u)\le f(x) \qquad
\mbox{in}\ \Omega, \ \mbox{in the viscosity sense}.
$$
\label{cor11new}
\end{cor}

\begin{rem} In the above theorem,
condition
 (\ref{ABC1}) is only needed to be satisfied,
for some $\bar r>0$,
in
$E_{\bar r}\setminus E$,
$E_{\bar r}=\{x\ |\
dist(x, E)<\bar r\}$, since
we can apply the theorem with
$\Omega=E_{\bar r}$.
\end{rem}

Our proof of Theorem \ref{thm11new} makes use of the following
maximum principle for functions satisfying (\ref{ABC1}).

\begin{prop}
For $n\ge 2$, $-1\le k\le n-2$,  let $E$ be a smooth closed
$k-$dimensional manifold in $\Bbb R^n$ and
$\Omega\subset \Bbb R^n$ be a domain.
Assume that
$u\in LSC(\overline \Omega\setminus E)$ satisfies
(\ref{ABC1}) and
$$
\inf_{\Omega\setminus E}u>-\infty.
$$
Then
$$
u\ge \inf_{ \partial \Omega \setminus E}u,\qquad \mbox{on}\ \Omega\setminus E.
$$
\label{prop1}
\end{prop}
Note that in the above, when $k=-1$, $E$ is understood as
$\emptyset$, the empty set;
 while for $k=0$, $E$ consists of finitely many points.

\begin{rem} 
The above proposition was proved in \cite{CLN1} 
under a stronger assumption that 
$u\in C^2(\Omega\setminus E)\cap C^0(\overline \Omega\setminus E)$.
The proof applies with minor modification.
\end{rem}

\noindent{\bf Proof of Theorem \ref{thm111new}.}\
Let $\varphi\in C^2(\Omega)$ satisfy
$$
\varphi \le u,\quad\mbox{in}\ \Omega,\qquad \mbox{and}\ \varphi(\bar
x)=u(\bar x),\ \mbox{for some}\
 \bar x\in \Omega.
$$
We only need to prove 
(\ref{aa1}).
If $\bar x\in \Omega\setminus E$, this follows from
(\ref{abcde}).  
 We may assume,
without  loss of generality, that
$\bar x=0\in E$.

For any fixed $\delta>0$, let
$$
\varphi^{ (\delta)}(x):= \varphi(x)-\frac \delta 2 |x|^2.
$$
Consider, for $0<\epsilon<1$,
$$
\varphi_\epsilon(x):= \varphi^{ (\delta)} +\epsilon d(x),
$$
where $d(x):=dist(x, E)$ denotes the distance function from $x$ to $E$.
Since $u$ is lowerconical on $E$,  
we have, for small $\epsilon>0$,  that
$$
\lambda(\epsilon):=-\inf_{\Omega\setminus E }(u-
\varphi_\epsilon)>0.
$$

 Since $\varphi_\epsilon=
\varphi^{ (\delta)}$ on $E$, and $u-
\varphi^{ (\delta)}>0$ on $\overline \Omega\setminus \{\bar x\}$,
we have, for small $\epsilon>0$,
 $$
u-\varphi_\epsilon>0\ \mbox{on}\ \partial \Omega, \qquad
 \liminf_{ x\in \Omega\setminus E, d(x)\to 0}
 (u-\varphi_\epsilon)(x)\ge 0.
 $$
Therefore,
 there
exists $x_\epsilon\in  \Omega\setminus E$ such that
$$
u-\tilde \varphi_\epsilon\ge 0\ \mbox{in}\
 \Omega\setminus E, \qquad
(u-\tilde \varphi_\epsilon)(x_\epsilon)=0,
$$
where
$$
\tilde \varphi_\epsilon(x):=
\varphi_\epsilon(x)-\lambda(\epsilon)
=\varphi^{ (\delta)}(x)+\epsilon d(x)-\lambda(\epsilon).
$$
Using the positivity of $u-\varphi^{ (\delta)}$ 
in $\overline \Omega \setminus \{0\}$,
we obtain from the above that
$$
\lambda(\epsilon)=-(u-\varphi_\epsilon)(x_\epsilon)=
\varphi^{ (\delta)}(x_\epsilon)
-u(x_\epsilon) +\epsilon d(x_\epsilon)\le \epsilon
d(x_\epsilon),
$$
and
$$
\lim_{\epsilon\to 0}
d(x_\epsilon)=\lim_{\epsilon\to 0}
|x_\epsilon|=0.
$$

By (\ref{abcde}), 
$$
F(x_\epsilon, \tilde \varphi_\epsilon(x_\epsilon),
\nabla \tilde \varphi_\epsilon(x_\epsilon),
\nabla^2 \tilde \varphi_\epsilon(x_\epsilon))\le f(x_\epsilon).
$$
We know from the above that, as $\epsilon\to 0$,
$$
x_\epsilon\to 0,\ \ \ 
 \tilde \varphi_\epsilon(x_\epsilon)\to  \tilde \varphi_\epsilon(0)
= \varphi(0),
\ \ \ 
 \nabla  \tilde \varphi_\epsilon(0) 
\to \nabla    \varphi(0).
$$
For
$$
\nabla^2 \tilde \varphi_\epsilon(x_\epsilon)=
\nabla^2 \varphi^{ (\delta)}(x_\epsilon)
+\epsilon \nabla^2 d(x_\epsilon),
$$
we use lemma 7.1 in 
\cite{CLN1} to obtain that $d(x)$ is pseudoconvex near $E$, i.e.
$$
\nabla^2 d(x_\epsilon)\ge O(1), \qquad
\mbox{as}\ \epsilon\to 0
$$
--- this is probably a known result.
Thus, using  (\ref{elliptic}), we have
\begin{eqnarray*}
f(0)&\ge&\lim_{\epsilon\to 0}f(x_\epsilon)
\ge \lim_{\epsilon\to 0} F(x_\epsilon, \tilde \varphi_\epsilon(x_\epsilon),
\nabla \tilde \varphi_\epsilon(x_\epsilon), 
\nabla^2 \tilde \varphi_\epsilon(x_\epsilon))\\
\\
&\ge & 
 \lim_{\epsilon\to 0} F(x_\epsilon, \tilde \varphi_\epsilon(x_\epsilon),
\nabla \tilde \varphi_\epsilon(x_\epsilon),
\nabla^2 \varphi^{ (\delta)}(x_\epsilon)
+O(\epsilon))\\
&=&F(0,  \varphi(0), \nabla
\varphi(0), \nabla^2 \varphi(0)+\delta I).
\end{eqnarray*}
Sending $\delta$ to $0$ in the above leads to the desired 
inequality  
(\ref{aa1}).
Theorem \ref{thm111new} is established.

\vskip 5pt
\hfill $\Box$
\vskip 5pt

Now we give the

\subsection{Proof of Theorem \ref{thm11new}}
It is clear that the extension $u$ is in
$LSC(\Omega)$.
We will prove that $u$
is lowerconical in $E$.

Fix any $\bar x\in E$ and any $\eta\in C^\infty(\Omega)$, 
we will show that (\ref{conical}) holds for every $\epsilon>0$.
We prove this by contradiction.  Suppose not, then
for some $\bar\epsilon>0$,
\begin{equation}
\bigg(  (u+\eta)(x)- (u+\eta)(\bar x) -\bar \epsilon
dist(x, E)\bigg)\ge 0,
\qquad \forall x\in \Omega.
\label{BBB1}
\end{equation}

 We may assume $\bar x=0\in E$, and
 the tangent space of $E$ at $0$ is spanned
by $e_{n-k+1}, \cdots, e_n$, where $e_1=(1, 0, \cdots, 0)$,
$\cdots, e_n=(0, \cdots, 0, 1)$ are the standard basis of
$\Bbb R^n$.  We write $x=(x', x'')$, where
$x'=(x_1, \cdots, x_{n-k})$ and
$x''=(x_{n-k+1}, \cdots, x_n)$.

For $x$ close to $0$, we have,  for some constant $C$,
\begin{equation}
\frac 34 |x'|-C|x''|^2\le
|x'|-C|x|^2\le
 d(x)\equiv dist(x, E)
\le |x'|+C|x|
\le
\frac 54 |x'|+C|x''|^2.
\label{48-1}
\end{equation}
The above fact
follows  easily from (70) in \cite{CLN1}.

It follows from (\ref{BBB1}) and (\ref{48-1}) that
 for some $\bar r,  \overline C>0$,
\begin{eqnarray*}
u(x)&\ge&
u(0)+\eta(0)-\eta(x)+\bar \epsilon
dist(x, E)=u(0)-\nabla \eta(0)\cdot x +O(|x|^2)
+\bar \epsilon
dist(x, E)\\
&\ge &
u(0)-\nabla \eta(0)\cdot x+
\frac {\bar  \epsilon} 2
|x'|-\overline C|x''|^2,\qquad x\in B_{\bar r}\setminus E.
\end{eqnarray*}

For $0<a<
\min\{ \bar r, \bar \epsilon/4\}$ which will be chosen later, let
$$
h(x):= u(0)-\nabla \eta(0)\cdot x+
\frac {\bar \epsilon}{ 4a}|x'|^2-(\overline C+1)|x''|^2+
\frac {a^2}2.
$$

On $\partial B_a(0)\setminus E$,
\begin{eqnarray*}
u(x)&\ge &u(0)-\nabla \eta(0)\cdot x+
  \frac {\bar \epsilon}{2a}|x'|^2
-\overline C|x''|^2 =h(x)+  \frac {\bar \epsilon}{4a}|x'|^2
+|x''|^2-\frac{a^2}2\ge h(x).
\end{eqnarray*}

By lemma 8.2  
in Appendix B of \cite{CLN1}
 and the assumption (\ref{ABC1}), there exists some
positive constant $\bar a>0$ such that if we further
require $0<a<\bar a$, we have
\begin{equation}
\sum_{i=1}^{k+2}\lambda_i(\nabla ^2 (u-h))
\le 0,\qquad
\mbox{in}\
\Omega\setminus E,\quad
\mbox{in the viscosity sense}.
\label{ddd1}
\end{equation}
Indeed,
let $\varphi\in C^2(\Omega\setminus E)$ satisfy
$$
\varphi \le u-h,\quad\mbox{in}\ \Omega\setminus E,
\qquad \mbox{and}\ \varphi(\bar
x)=(u-h)(\bar x),\ \mbox{for some}\
 \bar x\in \Omega.
$$
Then, by (\ref{ABC1}),
$$
\sum_{i=1}^{k+2}\lambda_i(\nabla ^2 h(\bar x)+
\nabla ^2 \varphi(\bar x))\le 0.
$$
Applying the above mentioned lemma in \cite{CLN1}, with 
 $l=k+2$, 
$$
D=  \frac {2a}{\bar \epsilon}
\nabla ^2 h(\bar x)=diag(1, \cdots, 1, -\delta_1, \cdots, -\delta_k),
\ \ \ 
\delta_1=\cdots=\delta_k= \frac{4a(\overline C+1)}{\bar \epsilon},
$$
$$
M=D+  \frac {2a}{\bar \epsilon}\nabla^2 \varphi(\bar x),
$$
we obtain
$$
\sum_{i=1}^{k+2}\lambda_i(\nabla ^2
\varphi(\bar x))
= \frac {\bar \epsilon}{ 2a}
\sum_{i=1}^{k+2}\lambda_i(M-D) 
\le \sum_{i=1}^{k+2}\lambda_i(M)
=  \frac {2a}{\bar \epsilon}
\sum_{i=1}^{k+2}\lambda_i(\nabla ^2 h(\bar x)+
\nabla ^2 \varphi(\bar x))\le 0.
$$
We have proved (\ref{ddd1}).
Thus, in view of proposition \ref{prop1} in \cite{CLN1},
$$
u-h\ge \inf_{ \partial B_a(0) \setminus E}(u-h)
\ge 0,
$$
and therefore
$$
\liminf_{x\to 0, x\in \Omega\setminus E}u(x)\ge h(0)=u(0)+\frac {a^2}2>u(0).
$$
A contradiction.
We have therefore proved (\ref{conical}).
Theorem \ref{thm11new} is established.

\vskip 5pt
\hfill $\Box$
\vskip 5pt

\section{Maximum principle}

In a bounded open
set $\Omega$ in $\Bbb R^n$ we consider two functions, $u, v$;
$u$ is in $LSC(\overline \Omega)$, and $v\in C^2(\Omega)\cap
C(\overline\Omega)
$.  The function $u$ is assumed to satisfy, in $\Omega$,
\begin{equation}
F(x,u,\nabla u, \nabla^2 u)\le
F(x,v,\nabla v, \nabla^2 v)\quad\mbox{in viscosity sense}.
\label{3.1}
\end{equation}
Here $F(x,s,p,M)$ is continuous and its derivatives
in $(s,p,M)$ are continuous. 
Concerning ellipticity of $F$ we assume 
here that $F$ may be degenerate elliptic, but
that there is a unit vector $\xi$, 
such that for all values of the arguments of $F$,
$$
F_{M_{ij}}\xi_i\xi_j>0.
$$
However we do not assume that this expression is
uniformly bounded by some positive constant.

\begin{thm}  \ (maximum principle)\  Assume
\begin{equation}
u\ge v\quad\mbox{on}\ \partial \Omega,
\qquad \inf_\Omega u>-\infty.
\label{3.2}
\end{equation}
Then $u\ge v$ in $\Omega$ provided
\begin{equation}
F_s(x,s,p,M)\le 0,\qquad\forall\ (s,x,p,M).
\label{3.3}
\end{equation}
\label{thm3.1}
\end{thm}

For a uniformly elliptic operator one knows that even if (\ref{3.3}) is not 
assumed, the conclusion $u\ge v$ holds if the volume
 of $\Omega$ is small.  However if there is
no uniform ellipticity this needs not hold.  Here is an

\begin{exam}
Let $\Omega=B_R$, $a=R^{-2}$,  and let
$$
F(x, u, \nabla u, \nabla^2 u):=
-e^{ -\Delta u}+1+u -\frac 12 x\cdot \nabla u.
$$
Then
$$
u(x):=-1+a|x|^2,  \ \ \mbox{and}\ v(x)\equiv 0
$$
satisfy
$$
F(x, u, \nabla u, \nabla^2 u)\le 0=F(x, v, \nabla v, \nabla^2 v),
\qquad \mbox{in}\ B_R,
$$
and
$$
u=v\ \mbox{on}\ \partial B_R.
$$
But
$$
u<v\ \mbox{in}\ B_R.
$$
  Note that $R$ may be arbitrarily small.
\label{example2.1}
\end{exam}

Before proving Theorem \ref{thm3.1} it is convenient to subtract $u$, and to 
consider $F(x,s,p,M)-F(x,v(x), \nabla v(x), \nabla^2 v(x))$
in place of $F$.  Then for the new $u$ and $F$ we have
\begin{equation}
F(x, u, \nabla u, \nabla^2 u)\le 0=
F(x,0,0,0)\quad
\mbox{in viscosity sense},
\label{3.1prime}
\end{equation}
and
\begin{equation}
u\ge 0\quad\mbox{on}\
\partial \Omega.
\label{3.2prime}
\end{equation}
From now on we assume $u$ satisfies 
(\ref{3.1prime}) and (\ref{3.2prime}).  Condition
(\ref{3.3}) continues to hold.

\noindent{\bf Proof of Theorem \ref{thm3.1}.}\
We may suppose that $F_{M_{11}}>0$.  We argue by contradiction.
Assume there is a point, which we take as
origin, where $u$ assumes its minimum value $-k$, $k>0$.

Suppose 
$$
\min_{\overline \Omega} x_1=-R.
$$

We use the comparison function
$$
h(x):= -k +k
e^{-\lambda (x_1+R) }-e^{ -\lambda R}
$$
with $\lambda>0$ to be chosen large.  We have 
$$
h(0)=-k, \qquad h<0\ \ \mbox{on}\ \partial \Omega.
$$

Move $h$ down, i.e. subtract a constant from $h$ so that it lies below $u$, then move it up,
to value 
$$
h-c_0, \qquad c_0\ge 0
$$
so that its graph first touches that of $u$ at some point $\bar x$.  Since 
$h<0$ on $\partial \Omega$, $\bar x\in \Omega$.  We have
$$
h_i= -k\lambda \delta_{i1} e^{ -\lambda(x_1+R) },\qquad
h_{ij} =k \lambda^2 \delta_{i1} \delta_{j1}  e^{ -\lambda(x_1+R) }.
$$
Because of (\ref{3.1prime}),
$$
F(\bar x, h(\bar x), \nabla h(\bar x), \nabla^2(\bar x))\le 0.
$$
Here all the arguments are bounded in absolute value 
by some constant independent of $\lambda$
and if we use the theorem of the mean, and the fact that $F(x, 0, 0, 0)=0$. we see that
\begin{equation}
0\ge a_{ij} h_{ij} +b_i h_i +ch =: I
\label{3.4}
\end{equation}
with $(a_{ij})$ uniformly positive definite, and all coefficients bounded in
absolute value.  In addition, by (\ref{3.3}),
\begin{equation}
c\le 0.
\label{3.5}
\end{equation}
Computing, we find
\begin{eqnarray*}
I&=& k\lambda^2 a_{11} e^{ -\lambda(x_1+R) } -k\lambda b_1
e^{ -\lambda(x_1+R) }+ c(h-c_0)\\
&\ge & k e^{ -\lambda(x_1+R) }[a_{11}\lambda^2
-\lambda b_1 +c];
\end{eqnarray*}
where we have used (\ref{3.5}).  But since $a_{11}>0$,
for large $\lambda$ this is positive, contradicts (\ref{3.4}).

\vskip 5pt
\hfill $\Box$
\vskip 5pt

\section{Strong maximum principle and Hopf Lemma for viscosity solutions}

We take up first the strong maximum principle.  Here
$$
u\ge v
$$
are functions defined in  $\Omega$,
an open and connected
subset of $\Bbb R^n$,
$u$ is in $LSC(\Omega)$ while $v$ is in $C^2(\Omega)$.
The function $u$ satisfies
$$
F(x,u,\nabla u,\nabla^2u)\le
F(x,v,\nabla v,\nabla^2v)\quad\mbox{in viscosity sense}.
$$

The nonlinear operator $F(x,s,p,M)$
 is continuous and of class $C^1$ in $(s,
p, M)$ for all values of the arguments.  $F$ is assumed to be elliptic,
i.e.
$$
(\frac {\partial F}{ \partial M_{ij} })\
\mbox{is positive definite},
$$
for all values of the arguments.  However $F$ is not assumed to
 be uniformly elliptic,
nor are $|F_s|$, $F_{p_i}|$ uniformly bounded.

\begin{thm} \ (Strong maximum principle)\
Let $u$ and $v$ be as above.  Suppose $u=v$ at some point in $\Omega$.  Then
$$
u\equiv v.
$$
\label{thm4.1}
\end{thm}

Before giving the proof, it is convenient
to change $u$ and $F$.  Namely, if we subtract $v$ from
$u$ and $F(x, v(x), \nabla v(x), \nabla^2 v(x))$ from $F$, we may then assume
 that
$$
u\ge 0
$$
and
\begin{equation}
F(x,u,\nabla u,\nabla^2u)\le 0
=F(x, 0, 0, 0)
\quad
\mbox{in viscosity sense}.
\label{4.2prime}
\end{equation}

From now on we assume that $u$ satisfies (\ref{4.2prime}).

\noindent{\bf
Proof of Theorem
 \ref{thm4.1}.}
\
We argue by contradiction.
Suppose the conclusion is false.  Since $u$ is LSC, and nonnegative,
the set where $u=0$ is closed.  Then there is an open ball $B$
of radius $R$,
with $\overline B\subset \Omega$,
with $u>0$ in $\overline B$ except
that $u(\hat x)=0$ at some point $\hat x$
on $\partial B$; we may suppose the center is the origin.

As in the classical proof of the strong maximum principle we make use
of a comparison function
\begin{equation}
h(x)=
E(x)-e^{ -\alpha R^2}, \quad
E(x):= e^{-\alpha |x|^2},
\qquad \alpha>0 \ \mbox{to be chosen large}.
\label{4.3}
\end{equation}
Then 
$$
h_i= -2\alpha x_i E, \qquad h_{ij}
=(4\alpha^2 x_ix_j -2\alpha \delta_{ij})E.
$$

Let
$A$ be an open ball, 
with $\overline A\subset \Omega$, having $\hat x$ as center
, of radius $\delta=\delta(\alpha)<R/2$ satisfying
\begin{equation}
\delta \alpha^{1/2}< \frac \pi {10}.
\label{4.6}
\end{equation}
Clearly
\begin{equation}
-1\le h\le 1, \quad |\nabla h|+|\nabla^2 h|\le C,
\qquad\mbox{in}\ A,
\label{3-1}
\end{equation}
where $C$ is some constant independent of $\alpha$.

Now, in the ball $\overline A$ we change $u$ and
$F$ to $\widetilde u$ and $\widetilde F$ to ensure that
$$
\widetilde F_{  \widetilde u }<0
$$
 for values of the arguments bounded, say, by $1$:
We set
\begin{equation}
u=\widetilde u\xi,
\qquad
\xi:= \frac 1\alpha \cos( (x_1-\hat x_1) \alpha^{1/2} )
\label{4.5}
\end{equation}
By (\ref{4.6}),
$$
 \cos( (x_1-\hat x_1) \alpha^{1/2} )
>\frac 12\quad\mbox{in}\ \overline A.
$$

Then set 
$$
\widetilde F(x, \widetilde u, \nabla  \widetilde u, \nabla^2  \widetilde u)
:= F(x,  \widetilde u \xi, \nabla (\widetilde u \xi),
\nabla^2  (\widetilde u \xi)),
$$
so that $\widetilde u$ satisfies 
\begin{equation}
\widetilde F(x, \widetilde u, \nabla  \widetilde u, \nabla^2  \widetilde u)
\le 0\quad
\mbox{in viscosity sense}.
\label{4.7}
\end{equation}

For some $\bar \epsilon=\bar \epsilon(\alpha)>0$, we have
$$
\widetilde u\ge \epsilon h,\qquad\mbox{on}\ \partial A,
\ \ \forall\ 0<\epsilon<\bar\epsilon.
$$
Now move $\epsilon h$ down, i.e. subtract a constant from it so that it lies below 
$\widetilde u$ in $\overline A$.  Then move it up, it becomes
$$
\epsilon h-c_0,\qquad 0\le c_0=c_0(\epsilon,\alpha)\le \epsilon,
$$
so that its graph first touches that of $\widetilde u$ at
some point $\bar x$.  Then, because of (\ref{4.7}), 
we have
$$
F(\bar x, (\epsilon h-c_0)\xi,
 \nabla ( (\epsilon h-c_0)\xi)(\bar x),
 \nabla ^2 ( (\epsilon h-c_0)\xi)(\bar x))\le 0.
$$
Because of (\ref{3-1}), we can fix some constant
$\epsilon\in (0, \bar\epsilon(\alpha))$
 such that
$$
\max\{ |(\epsilon h-c_0)(\bar x)|,
|\nabla (\epsilon h-c_0)(\bar x)|,
|\nabla^2(\epsilon h-c_0)(\bar x)|\}\le 1.
$$
Thus, at $(\bar x, (\epsilon h-c_0)(\bar x),
 \nabla  (\epsilon h-c_0)(\bar x),
 \nabla ^2  (\epsilon h-c_0)(\bar x))$,
$F_{M_{ij}}$ is uniformly
positive definite and $|F_{p_i}|$, $|F_s|$ are
bounded independent of $\alpha$.
By the theorem of the mean we find that at $\bar x$,
\begin{equation}
F_{ M_{ij} } (\epsilon h-c_0)_{ij} 
+[2F_{ M_{ij} } \xi_j\xi^{-1} + F_{p_i}]
(\epsilon h-c_0)_i
+c (\epsilon h-c_0)\le 0.
\label{4.8}
\end{equation}

\noindent {\bf Claim.}\ $c<0$.

\noindent {\bf Proof.}\ Here
$$
 c\xi =F_{ M_{ij} }\xi_{ij} +F_{p_i}\xi_i+F_s\xi,
$$
where the arguments in $F$ and its derivatives are all bounded, independent of $\alpha$.  Also $F_{M_{ij}}$ is uniformly 
positive definite and $|F_{p_i}|$, $|F_s|$ are 
bounded independent of $\alpha$.  Hence, for large
$\alpha$,
$$
 c\xi = -F_{M_{11}}\cos[(x_1-\hat x_1)\alpha^{1/2}]
+O(\alpha^{-1/2})<0.
$$

Since $c<0$ we see that
$$
c(\epsilon h-c_0)= c\epsilon E-c(\epsilon e^{-\alpha k^2}+c_0)
>c \epsilon E .
$$
Inserting this in (\ref{4.8}) we infer that
$$
F_{ M_{ij} }  h_{ij}
+ [2F_{ M_{ij} } \xi_j\xi^{-1} + F_{p_i}]
h_i
+c  E <0,
$$
i.e.
$$
 F_{ M_{ij} } (4\alpha^2 x_ix_j-2\alpha \delta_{ij})E
+[2F_{ M_{i1} } \xi_1\xi^{-1} + F_{p_i}]
(-2\alpha x_i E)
+c E<0.
$$

Since  $|x|$ is bounded
away from zero in $\overline A$, we see from the above that
 for some positive constant
$c_1, c_2$ independent of $\alpha$,
$$
(c_1 \alpha^2 -c_2\alpha -c_2 \alpha ^{3/2}-c_2\alpha -c_2)E<0.
$$
But this is impossible for $\alpha$ large.

\vskip 5pt
\hfill $\Box$
\vskip 5pt

Next

\begin{thm}\ (Hopf Lemma)\ Let $u$ and $v$ be as above, with $u>v$
in $\Omega$, and suppose 
$$
u(\hat x)=v(\hat x)
$$
at a boundary point $\hat x$ near which $\partial \Omega$ is $C^2$.
Then, if $\nu$ is the unit interior normal to
$\partial \Omega$ at $\hat x$,
$$
\liminf_{s\to 0^+}
\frac { (u-v)(\hat x+s\nu)}s>0.
$$
\label{thm4.2}
\end{thm}

\noindent{\bf Proof.}\  As before, by considering $u-v$ in place of $u$, etc. we may suppose
$$
u>0\quad\mbox{in}\ \Omega,\qquad u(\hat x)=0,
$$
$$
F(x, u, \nabla u, \nabla^2 u)\le 0
\quad
\mbox{in viscosity sense
in}\ \Omega.
$$

Let $B$ be a ball of radius $R$, in $\Omega$ with $\hat x$ on its boundary.  We take the
origin as center of $B$.  
 We use the same comparison function $h$ as in
(\ref{4.3}).   
As in the proof of Theorem \ref{thm4.1}, let
$A$ be an open  ball with center at $\hat x$, with radius $\delta(\alpha)<R/2$
satisfying
(\ref{4.6}).
We work in the region
$$
D=B\cap  A.
$$
In $\overline B$ we introduce as before the function
$\widetilde u$ defined as in (\ref{4.5}), and
$\widetilde F$.
 
Again, for small $\epsilon>0$ we have
\begin{equation}
\widetilde u>\epsilon h\quad \mbox{on}\ \partial D.
\label{4.11}
\end{equation}
Now, argue as in the proof of Theorem \ref{thm4.1}.
Move $\epsilon h$ down and then up so that it becomes
$$
\epsilon h-c_0, \qquad 0\le c_0=c_0(\epsilon, \alpha)\le \epsilon,
$$
and so that its graph first touches that of
$\widetilde u$ from below
at some point $\bar x$.

\noindent{\bf Claim.}\ $c_0=0$.

If that is the case, then $\widetilde u\ge \epsilon h$ in
$\overline D$ and the desired conclusion,
$$
\liminf_{s\to 0^+}
\frac { u(\hat x+s\nu)}s>0,
$$
follows.

\noindent{\bf Proof of Claim.}\ Suppose not, suppose $c_0>0$.  Because  of (\ref{4.11}), $\bar x$ is in $D$.  Then arguing exactly as in the proof of Theorem 
\ref{thm4.1} we are led to a contradiction.

\vskip 5pt
\hfill $\Box$
\vskip 5pt

At the end of this section we point out that 
Theorem \ref{thm3.1} can be deduced from Theorem \ref{thm4.1} as follows.

\noindent{\bf Theorem \ref{thm3.1} as a consequence of
Theorem \ref{thm4.1}.}\
Suppose the contrary, then
\begin{equation}
\inf_\Omega (u-v)<0.
\label{H1-1}
\end{equation}
Move $v$ down so that its graph lies below that of $u$, then move it up to
value 
$$
v-c_0, \quad c_0\ge 0,
$$
so that its graph first touches that
of $u$ at some point $\bar x\in \overline \Omega$.  Namely,
$$
u\ge v-c_0\quad\mbox{in}\ \Omega,
\qquad u(\bar x)=v(\bar x)-c_0.
$$
By (\ref{3.2}) and (\ref{H1-1}),
$c_0>0$ and $\bar x\in \Omega$.
By (\ref{3.3}), we have
$$
F(x,v,\nabla v, \nabla ^2 v)\le F(x,v-c_0, \nabla (v-c_0),
\nabla^2 (v-c_0)).
$$
Thus, in view of (\ref{3.1}),
$$
F(x,u,\nabla u, \nabla^2 u)\le
F(x,v-c_0, \nabla (v-c_0),
\nabla^2 (v-c_0)) \quad\mbox{in}\ 
\Omega,\ \mbox{in viscosity sense}.
$$
Applying Theorem \ref{thm4.1} to $u$ and $v-c_0$, we infer that
$u= v-c_0$ in the connected component of $\Omega$ containing
$\bar x$.  This violates $u\ge v$ on $\partial \Omega$.

\vskip 5pt
\hfill $\Box$
\vskip 5pt

\section{ Strong maximum principle and Hopf Lemma for viscosity
solution of nonlinear parabolic equation}

In this section we extend the strong maximum principle,
Theorem \ref{thm4.1}, to nonlinear 
parabolic operators.  

In the closure $\overline \Omega$ of a domain $\Omega$ in
$\Bbb R^{n+1}$, $(x,t)$ space,
$x\in \Bbb R^n$, $t\in \Bbb R$, we consider
two functions
$$
u\ge v,
$$
$u$ is lower semicontinuous (LSC) while $v\in C^2$;
$u$ satisfies, in the viscosity sense, in $\Omega$,
\begin{equation}
F(x,t,u,\nabla u, \nabla^2 u)-u_t\le
F(x,t,v,\nabla v, \nabla^2 v)-v_t.
\label{5.1}
\end{equation}

Here $\nabla$ and $\nabla^2$ represent first and second derivatives with
respect to the $x-$variables.  $F(x,t,s,p,M)$
is as in
section 3: $F$ is continuous and of class $C^1$ in $(s,p,M)$ for all values
of the arguments.  $F$ is assumed to be elliptic, i.e.
$$
(\frac{\partial F}{ \partial M_{ij}})
\ \ \mbox{is positive definite},
$$
for all values of the arguments.
However $F$ is not assumed to be uniformly elliptic,
nor are $|F_s|$, $|F_{p_i}|$ uniformly bounded.

\noindent {\bf Setup.}\ We assume that $\Omega$ lies in $
\{t<T\}$ for some $T$ and that $\partial \Omega$ includes a relatively open
subset $\Sigma$ on the hypersurface $\{t=T\}$.

For every point $P=(x_0, t_0)
\in \Omega\cup \Sigma$, we
denote by $C_P$ the arcwise connected component,
containing $P$, of points $(x, t_0)$, in $\Omega\cup \Sigma$.
We emphasize that $C_P\subset \{t=t_0\}$.

We also denote by $S_P$ the set of points in $\Omega$
which may be connected to $P$ by a continuous 
curve on which the $t-$coordinate is nondecreasing.

We require also that 
\begin{equation}
\mbox{
(\ref{5.1}) holds not only on
}\ \Omega, \ \mbox{
but also at points of
}
\ \Sigma.
\label{5.1prime}
\end{equation}

\begin{rem}  This makes sense.  It would not make sense if we require (\ref{5.1}) to hold at point $(\bar x, \bar t)$ on a lower boundary 
point of $\Omega$, i.e. when $\Omega$
lies in $\{t>\bar t\}$.
\end{rem}

The main result of this section, the
parabolic strong maximum principle, is
\begin{thm}
Let $\Omega$, and $u\ge v$ be as above.  If 
$u(P)=v(P)$ at a point on
$\Sigma$ then
$$
u\equiv v\quad \mbox{in}\ C_P\cup S_P.
$$
\label{thm5.1}
\end{thm}

Before starting the proof it is convenient, as in the elliptic case,
to change $u$ and $F$.  If we subtract $v$ from $u$, and
$F(x, t, v, \nabla v, \nabla^2 v)$ from
$F$ then we may suppose 
\begin{equation}
u\ge 0\equiv v,\qquad \mbox{in}\ \Omega\cup \Sigma,
\label{5.3a}
\end{equation}
\begin{equation}
F(x,t,u,\nabla u, \nabla^2 u)-u_t\le 0
=F(x,t,0,0,0),
\qquad \mbox{in}\ \Omega\cup \Sigma,
\ \mbox{in the viscosity sense.}
\label{5.3b}
\end{equation}

From now on we assume $u$ satisfies (\ref{5.3a}) and (\ref{5.3b}).

The proof follows that of the parabolic strong maximum principle
in Nirenberg \cite{N}, with modifications for viscosity supersolution,  It makes use of
several lemmas.

\begin{lem} Consider $u$ satisfying
(\ref{5.3a}) and (\ref{5.3b}); suppose
also $u>0$ in a ball,
$B$, with $\overline B\subset \Omega$,
$$
u=0\ \ \mbox{at a point}\ 
P\ \mbox{on}\ \partial B.
$$
Then necessarily, the vector
$(0,\cdots, 0,1)$ is normal to $\partial B$ at $P$.
\label{lem5.1}
\end{lem}

We postpone the proof; first make some corollaries:
\begin{cor}
If $u>0$ in a subdomain $G$ of $\Omega$, with
$\overline G\subset \Omega$ and $u(P)=0$ at
a point $P$ on $\partial G$ where $\partial G$ is smooth, then
$(0,\cdots, 0,1)$ is normal to 
 $\partial G$ at $P$.
\label{cor5.1}
\end{cor}

This follows from Lemma \ref{lem5.1} by just taking a ball $B$ in $G$ with $P$ on its boundary.

\begin{cor}
Let $u$ satisfy (\ref{5.3a}) and (\ref{5.3b}).  If $u(P)=0$ for
some point $P\in \Omega$ then
$$
u\equiv 0\ \ \mbox{on}\ C_P.
$$
\label{cor5.2}
\end{cor}

\noindent{\bf Proof.}\
Suppose not, suppose $u(Q)>0$ for some
$Q\in C_P$.  Join $Q$ to $P$ by a continuous curve
on $C_P$.  As we traverse
 the curve from $Q$ to $P$, let $\bar P$ be the first point where $u=0$; it maybe $P$.  Let $\bar Q$ be a point on the curve so close to $\bar P$ that
$B_{ 8 |\bar Q-\bar P|}(\bar Q)\subset \Omega$.

Since $u$ is LSC there is a small vertical segment,
i.e. parallel to $t-$axis, of length $2\epsilon$,
$0<\epsilon<|\bar Q-\bar P|$,
with center at $\bar Q$, where $u>0$.  If $\bar Q
=(\bar x, \bar t)$, the closed ellipsoid
$$
E_a:= \{(x,t)\ |
\ (t-\bar t)^2+ a^{-2}|x-\bar x|^2
\le \epsilon^2\}
$$
lies in $\Omega$, provided
 $0<a\epsilon< |\bar Q-\bar P|$.

For small $a>0$, $u>0$ in $E_a$.  
Now increase $a$, as we do so, $u$ remains positive
in $E_a$.  This follows from Corollary \ref{cor5.1}.
Finally, $u>0$ in $E_a$ for $a=\bar a
:=|\bar Q-\bar P|/\epsilon$.
But $\bar P$ lies on the boundary of $E_{\bar a}$.  Contradiction.

\vskip 5pt
\hfill $\Box$
\vskip 5pt

We will also use
\begin{lem}
Let $u$ satisfy (\ref{5.3a}) and(\ref{5.3b}).  Suppose
$u>0$ in open half ball $D$ in
$t<T_0\le T$ centered 
at $(\bar x, T_0)$.  Then $u>0$ also 
on the relatively open part of the flat boundary
of $D$, where $t=T_0$.
\label{lem5.2}
\end{lem}

Before proving the lemmas we first show 
how they give the

\noindent{\bf Proof of Theorem \ref{thm5.1}.}\
Suppose it does not hold, i.e. there is
some point $Q$ in $C_P\cup S_P$ where $u(Q)>0$.  Without loss of
generality, since $u$ is LSC, we may suppose that
$Q$ is in $S_P$.  Join $Q$ to $P$ by a continuous 
curve $\Gamma$ on which $t$ is nondecreasing.  As we traverse 
the curve from $Q$, let $P_0
=(x_0, t_0)$ be the first point where $u=0$ (it may be $P$).  Let
$D$ be an open half ball with center at $P_0$
whose closure lies in $(\Sigma\cup \Omega)\cap\{t\le t_0\}$.
Near the end of $\Gamma$ it lies in $D$. 
By Corollary \ref{cor5.2}, $u>0$ in $D$. Then, by Lemma \ref{lem5.2},
$u(P_0)>0$.  Contradiction.  Theorem \ref{thm5.1}
is proved.

\vskip 5pt
\hfill $\Box$
\vskip 5pt

We now prove the lemmas.  First

\noindent{\bf Proof of Lemma \ref{lem5.2}.}\
Suppose the conclusion is false, we suppose $u(P) 
=0$ at a point $P=(x_0, t_0)$ on the relatively open part of the flat 
boundary of $D$.  We may take $x_0$ to be the origin.  In addition, for
convenience, we take $T_0=0$.  So $P=\{0\}$.

Near the origin we introduce the comparison function
$$
h=-\alpha t-|x|^2,\quad
0<\alpha\ \mbox{to be chosen large}.
$$
In the cutoff region
$$
D_\alpha:=\{(x,t)\
|\ -\alpha^{-2}<t<0,\
\alpha t+|x|^2<0\},
$$
we have
\begin{equation}
|x|<\frac 1{\sqrt{\alpha}},\ \ 0<h<\frac 1\alpha,\ \ 
h_t=-\alpha, \ \ 
h_i=-2x_i,\ \ h_{ij}=
-2\delta_{ij}.
\label{5.4}
\end{equation}

Since $u>0$ on the boundary of $D_\alpha$ where
also $t=-\alpha^{-2}$,
$u>\epsilon h$ there for small $\epsilon>0$
(the smallness of $\epsilon $ may depend on $\alpha$).  Thus,
on $\partial D_\alpha$, $u\ge \epsilon h$,
with equality only at $\{0\}$.

Now, move $h$ down, i.e., subtract a positive constant
from $h$, so that it lies below $\widetilde u$ in
$\overline D_\delta$.  Then move it up to
$\epsilon h-c_0$, $0\le c_0\le \epsilon/\alpha$, so that 
its graph first touches that of $u$ at some point 
$(\bar x, \bar t)\in D_\alpha\cup \{0\}$.

Since $u$ satisfies (\ref{5.1prime}),
at $(\bar x, \bar t)$,
$$
F(x,t,\epsilon h-c_0, \epsilon \nabla h,
 \epsilon \nabla^2 h)-\epsilon h_t\le 0.
$$
We will show that this cannot hold for $\alpha$
large.  For $0<\epsilon$ small, the arguments $(\epsilon h-c_0,
\epsilon \nabla h, \epsilon \nabla^2 h)$
 are all bounded in absolute value by $1$.  Hence, by
the ellipticity of $F$ and the fact that
$F(x,t,0,0,0)=0$,
$$
0\ge -\epsilon h_t +\epsilon a_{ij}
h_{ij}+\epsilon b_i h_i+c(\epsilon h-c_0)=: J,
$$
with $(a_{ij})$ uniformly 
positive definite, and all coefficients bounded in absolute value.

We now compute $J$. Using (\ref{5.4}) and the inequalities
$$
h\le \frac 1\alpha,\qquad 0\le c_0\le \epsilon/\alpha,
$$
we find for a fixed constant $C$ independent of $\alpha$ that
$$
\frac J\epsilon \ge \alpha -C-\frac C\alpha>0\quad\mbox{for}\ \alpha\ \mbox{large}.
$$
Contradiction.

\vskip 5pt
\hfill $\Box$
\vskip 5pt

Now, 

\noindent{\bf Proof of Lemma \ref{lem5.1}.}\
By taking a smaller ball inside $B$ 
with $P$ on its boundary, we may suppose that
$$
u>0\ \mbox{on}\ \overline B\ \mbox{except at}\ P.
$$
As usual we argue by contradiction.  Suppose the conclusion is false.  
We may suppose that the origin is the center of $B$,
 its radius is $R$, and $P=(\hat x, \hat t)$, $\hat x\ne 0$.
We use the comparison function
$$
h= E-e^{-\alpha R^2},\qquad E=e^{-\alpha (|x|^2+t^2)  }.
$$
We have
$$
h_t=-2\alpha t E, \ \ 
h_i=-2\alpha x_i E, \ \
h_{ij}=(-2\alpha \delta_{ij}+4
\alpha^2 x_ix_j)E,
$$

Let $A$ be a small ball centered at $P$, with 
radius $\delta<|\hat x|/2$ and $\overline A\subset \Omega$.
We require $\delta=\delta(\alpha)$ to be small, namely we require
$$
\delta \alpha^{1/2}<\frac \pi{10},
$$
so that
 $$
\cos[(x_1-\hat  x_1)\alpha^{1/2} ]>
\frac 12\qquad \mbox{in}\ \overline A.
$$

Now, in the ball $\overline A$ we change $u$ and $F$ to $\widetilde u$ and $\widetilde
F$, to ensure that
$$
\widetilde F_{ \widetilde u}<0
$$
for values of the arguments bounded, say,
by $1$.  Namely,
 we set
$$
u=\widetilde u\xi,
\quad \xi=\frac {  \cos[(x_1-\hat x_1)\alpha^{1/2} ]}\alpha.
$$

Then we set
$$
\widetilde F(x,t,\tilde u, \nabla \widetilde u,  \nabla^2 \widetilde u)
=\frac 1\xi F(x,t,  \widetilde u
\xi, \nabla ( \widetilde u\xi), \nabla^2 
( \widetilde u\xi) ),
$$
so $ \widetilde u$ satisfies
$$
\widetilde F(x,t,\tilde u, \nabla \widetilde u,  \nabla^2 \widetilde u)
-\tilde u_t\le 0
\qquad
\mbox{in viscosity sense}.
$$

For some $\bar\epsilon=\bar\epsilon(\alpha)>0$,
 we have
$$
\widetilde u\ge \epsilon h\quad\mbox{on}\
\partial A,\quad\forall\
0< \epsilon< \bar\epsilon.
$$
As we did in the proof of Lemma \ref{lem5.2}, move
$\epsilon h$ down so that it lies below $\widetilde
u$, and then move it up, so that it becomes
$$
\epsilon h-c_0,\qquad 0\le c_0
=c_0(\epsilon, \alpha)\le \epsilon,
$$
so that its graph first touches that of
$\widetilde u$ at some point $(\bar x, \bar t)$.  Then at
$(\bar x, \bar t)$,
$$
I:= F(\bar x, \bar t,
(\epsilon h-c_0)\xi, \nabla ((\epsilon h-c_0)\xi),
\nabla^2 ((\epsilon h-c_0))\xi))
-\epsilon \xi h_t\le 0.
$$
By the theorem of the mean we find that,
at $(\bar x, \bar t)$,
$$
F_{ M_{ij} } (\epsilon h-c_0)_{ij}
+[2F_{ M_{ij} } \xi_j\xi^{-1} + F_{p_i}]
(\epsilon h-c_0)_i
+c (\epsilon h-c_0)-\epsilon h_t\le 0,
$$
where 
 $$
c\xi =F_{ M_{ij} }\xi_{ij} +F_{p_i}\xi_i+F_s\xi,
$$
and
 the arguments in $F$ and its derivatives are
all bounded independent of $\alpha$; also,
$F_{ u_{ij}}$ is uniformly positive definite and
$|F_{p_i}|$, $|F_s|$ are bounded.

We claim that $c<0$.

For large $\alpha$,
$$
c\xi =-F_{ M_{11}}\cos[ (x_1-\hat x_1)\alpha^{1/2}]
+O(\alpha^{-1/2})<0,\qquad
\mbox{for}\ \alpha\ \mbox{large}.
$$
Using the fact that $c<0$, we
argue as in the proof of Theorem
 \ref{thm4.1} to obtain
$$
 F_{ M_{ij} } (4\alpha^2 x_ix_j-2\alpha \delta_{ij})E
+[2F_{ M_{i1} } \xi_1\xi^{-1} + F_{p_i}]
(-2\alpha x_i E)
+c E-h_t<0.
$$
Since  $|x|$ is bounded
away from zero in $\overline A$, we see from the above that
 for some positive constant
$c_1, c_2$ independent of $\alpha$,
$$
(c_1 \alpha^2 -c_2\alpha -c_2 \alpha ^{3/2}-c_2\alpha -c_2-c_2\alpha)E<0.
$$
But this is impossible for $\alpha$ large.
Contradiction.  Lemma \ref{lem5.1} is proved.

\vskip 5pt
\hfill $\Box$
\vskip 5pt

Using similar arguments we now prove a parabolic Hopf Lemma for viscosity supersolutions.

Consider $\Omega$ and $\Sigma$, and $u, v$ as above, with 
$$u>v\quad \mbox{in}\ \  \Omega\cup \Sigma.
$$
We will prove the parabolic Hopf Lemma at a point, which we take to be the 
origin $\{0\}$, on
$\partial \Sigma$.  

$\overline \Omega\setminus \Sigma$ is called the parabolic boundary,
$P\partial\Omega$ of $\Omega$, and we assume that it is of class $C^2$ near
$\{0\}$.  For convenience we suppose
that
$\nu=(0, \cdots, 0, 1, 0)$ is the inner normal 
to $\partial \Sigma$ (of class $C^2$) at $(0,0)$, and we denote $x_n$ by $y$.
Sometimes we use $(x,y,t)$ with $x=(x_1,
\cdots, x_{n-1})$.  We assume that the interior
normal to $P\partial\Omega$ at $\{0\}$ is not $(0, \cdots, 0, 1)$.

We now suppose 
$
u>v$ in $\Omega\cup \Sigma$ and $u(0,0)=v(0,0)$, and we assume 
(\ref{5.1prime}), i.e.
$$
F(x,t,u,\nabla u, \nabla^2 u)-u_t
\le F(x,t,v,\nabla v,\nabla^2 v)
-v_t\quad \mbox{in}\ \Omega\cup \Sigma
\quad\mbox{in viscosity sense}.
$$

\begin{thm}  \ (Parabolic Hopf Lemma)\
Under the conditions above,
\begin{equation}
\liminf_{ s\to 0^+} \frac{ (u-v)(s\nu) 
} s>0.
\label{5.13}
\end{equation}
\label{thm5.2}
\end{thm}

\begin{rem} It will be clear from the proof that (\ref{5.13}) will also hold 
for any unit vector $\nu=(\nu_1, \cdots,
\nu_{n+1})$ at $\{0\}$ which 
points into $\Omega\cup \Sigma$ and is not tangent 
to $P\partial \Omega$, so, $\nu_{n+1}\le 0$. 
\end{rem}

As before, by considering $u-v$ in place of $u$, and subtracting $F(x,t,v,
\nabla v, \nabla^2 v)$ from $F$ we may assume $v\equiv 0$ and
$$
F(x,t,u,\nabla u, \nabla^2 u)-u_t\le 0=F(x,t, 0, 0, 0)
\quad\mbox{in}\ \Omega\cup \Sigma, \quad \mbox{in viscosity sense}.
$$

\noindent{\bf Proof of Theorem \ref{thm5.2}.}\ By restricting $\Omega$ we may assume that near $\{0\}$,
$\partial\Sigma$ is given by 
$$
y=d |x|^2,\qquad d>0,
$$
and that for some constant $b>0$, the domain
$$
\widehat \Omega=
\{(x,y,t)\ |\ t<0,\
y>d|x|^2-bt\},
$$
near the origin, lies in $\Omega$.  By decreasing $d$ and
increasing $b$ we may suppose that for the resulting $\widehat \Omega$, which we now call $\Omega$,
\begin{equation}
u>0\ \ \mbox{on}\ P\partial\Omega\ \mbox{except at}\ \{0\}.
\label{5.15}
\end{equation}

We will take $b$ to be large.

Next we introduce the comparison function
$$
h= y-d|x|^2 -bt.
$$
With 
$$A=\mbox{ball centered at origin with radius}\
\delta
\ \mbox{small},
$$
we consider $u$ and $h$ in the region
$$
G=\Omega\cap A.
$$
Since (\ref{5.15}) holds we see that for some $0<\epsilon$ small,
$$
u>\epsilon h\ \ \mbox{on}\ P\partial G.
$$
The desired conclusion (\ref{5.13}),
$$
\liminf_{ s\to 0^+} \frac{u(s\nu)
} s>0,
$$
will follow if we can show that 
$u\ge \epsilon h$ on $\overline G$.

To achieve this we argue as before: lower $\epsilon h$ so that it lies below $u$ in
$\overline G$ and then raise it to
$$
\epsilon h-c_0,
$$
until its graph first touches that of $u$.  We claim
this must happen for $c_0=0$, which would prove 
$$
u\ge \epsilon h.
$$
Suppose not, suppose $c_0>0$ and that the point of contact is
$(\bar x, \bar t)$.  Clearly
$(\bar x, \bar t)$ is not on $P\partial\Omega$; $\bar t$ might be zero.  At
$(\bar x, \bar t)$ we have
$$
F(x,t, h, \nabla h, \nabla^2 h)-\epsilon h_t\le 0.
$$
All arguments in $F$ are bounded by $1$, for $\epsilon$ small, so that we may infer, as before, that
$$
0\ge a_{ij}h_{ij}+b_ih_i+ch-h_t.
$$
With the operator on the right uniformly elliptic
and with coefficients uniformly
bounded.  Thus for some $C$ independent of $b$,
$$
0\ge -C +b >0\quad\mbox{for}\ b\ \mbox{large}.
$$
Contradiction.  Then $c_0=0$, i.e.
$$
 u\ge \epsilon h.
$$

\vskip 5pt
\hfill $\Box$
\vskip 5pt

\section{A strengthened Hopf Lemma 
for viscosity solution of
 parabolic equations}

In this section we extend Lemma \ref{lemA-0new} to parabolic equations.
The result is not used in this paper.  On the other hand it
is
useful when extending Theorem \ref{thmnew1} to parabolic equations.
We plan to extend Theorem \ref{lem1}-\ref{thm11new} to  parabolic equations
in a forthcoming paper.

Let $\Omega\subset \Bbb R^n$ be a domain with $C^2$ boundary,
$0<T<\infty$.  Assume that 
$(a_{ij}(x, t))$,  $b_i(x, t)$ and $c(x, t)$
are  functions in  $L^\infty(\Omega\times (0, T])$  satisfying,
for some positive constants $\lambda$ and $\Lambda$,
\begin{equation}
|a_{ij}(x,t)|+|b_i(x, t)| +|c(x, t)|\le \Lambda,
\ \
a_{ij}(x, t)\xi_i\xi_j\ge \lambda|\xi|^2,\ \ \ \forall\
x\in \Omega, \ 0<t<T,\ \xi\in \Bbb R^n.
\label{aij-parabolic}
\end{equation}
We will use the notation
$$
Lu:= a_{ij}(x, t)\partial_{ij}u+b_i(x,t)\partial_iu+c(x,t)u.
$$

\begin{thm} For $0<T_1<T<\infty$ and $0<\delta<1$,
let
$(a_{ij}(x, t))$,  $b_i(x, t)$ and $c(x, t)$
be  $L^\infty(\Omega\times (0, T])$   functions satisfying
(\ref{aij-parabolic}) with $\Omega=B_1$
for some positive constants $\lambda$ and $\Lambda$.
There exist
some positive constants $\epsilon,  \mu>0$
which depend only on $n, \lambda, \Lambda, \delta, T_1, T$, such that
if  
  $u\in
LSC(\overline B_1\times (0, T])$ satisfies
\begin{equation}
(L-\partial_t)
 u\le   \epsilon,
\quad \mbox{in}\ B_1\times (0, T],
 \
\mbox{in viscosity sense},
\label{A1-parabolic}
\end{equation}
\begin{equation}
u(x,0)\ge 1,\qquad \mbox{for}\ |x|\le \delta,
\label{A2-parabolic}
\end{equation}
\begin{equation}
u\ge 0,
\quad \mbox{on}\ P\partial(B_1\times (0, T]).
\label{A3-parabolic}
\end{equation}
Then
\begin{equation}
u(x,t)\ge
  \mu (1-|x|),\qquad  \mbox{on}\ B_1
\times [T_1, T].
\label{A4-parabolic}
\end{equation}
\label{maximum-parabolic}
\end{thm} 
Recall that $P\partial(B_1\times (0, T])$ denotes the 
parabolic boundary of $B_1\times (0, T]$, i.e.
$$
P\partial(B_1\times (0, T])=(\overline B_1\times\{0\})
\cup( \partial B_1 \times [0, T]).
$$

\noindent{  \underline{Note}.}\
 The function $u$  may actually be negative somewhere.

\medskip

\begin{thm} Let $\Omega$ be a domain of $\Bbb R^n$ with
$C^2$ boundary, and let $B\subset \Omega$ be a ball.  
For $0<T_1<T<\infty$,
let
$(a_{ij}(x, t))$,  $b_i(x, t)$ and $c(x, t)$
be  $L^\infty(\Omega\times (0, T])$   functions satisfying
(\ref{aij-parabolic}) for some positive constants $\lambda$ and $\Lambda$.
There exist
some positive constants $\epsilon,  \mu>0$
which depend only on $n, \lambda, \Lambda, \Omega,$
the radius of $B$, $T_1, T$, such that
if
  $u\in
LSC(\overline \Omega\times (0, T])$ satisfies
$$
(L-\partial_t)
 u\le   \epsilon,
\quad \mbox{in}\ \Omega \times (0, T],
 \
\mbox{in viscosity sense},
$$
$$
u(x,0)\ge 1\quad\mbox{for all}\ x\in B,
$$
\begin{equation}
u\ge 0\quad\mbox{on}\ P\partial(\Omega\times (0, T]). 
\label{ABC2}
\end{equation}
Then
$$
u(x, t)\ge
 \mu dist(x, \partial \Omega), \qquad  \mbox{on}\ \Omega\times [T_1, T].
$$
\label{thm5.2new}
\end{thm}

\noindent{\bf Proof of Theorem \ref{maximum-parabolic}.}\ 
We only need to prove that there exists some
constant $\overline T$ depending
only on $n, \lambda, \Lambda, \delta$ such that 
the theorem holds under an additional 
assumption that $T\le \overline  T$.
Indeed, for general $T$, we  fix a positive integer
$m$ so that $\frac Tm<\overline T$, and then apply
the result on $[0, \frac Tm]$, $[\frac Tm,
\frac {2T}m]$, ..., $[\frac {(m-1)T}{m}, T]$ successively.
We leave the simple details to readers.

In the following we will assume 
that $T\le  \overline T$, and we will determine the value of
$\overline T$ later.

We may assume without loss of generality that
$c(x,t)\le 0$ for all $|x|<1$ and $0<t<T$.
This can be achieved by working with
$$
\tilde u(x,t)= e^{ -2\Lambda t} u(x,t),
$$
since (\ref{A1-parabolic}) implies
$$
(\widetilde L-\partial_t)\tilde u:=
(L-2\Lambda-\partial_t)\tilde u  =   
e^{ -2\Lambda t} (L -\partial_t) u\le\epsilon,
$$
and $\widetilde L$ has
$\tilde c=c-2\Lambda\le -\Lambda< 0$.

Consider the comparison function
$$
h(x,t):=\frac 1 D(E-F),
\ \ 
E= (t+a)^{-k} 
e^{ -\frac {\alpha |x|^2}{ t+a}  },
\ \
F=(T+a)^{-k}
e^{ -\frac {\alpha }{ T+a}  },
$$
with
$$
\alpha =\frac 1\lambda,
$$
and 
$$
\frac {T+a}{ 4a}= \frac 2{ \delta^2}.
$$
Thus
\begin{equation}
a:=(8\delta^{-2}-1)^{-1} T.
\label{7}
\end{equation}
Next we require that
\begin{equation}
h(\frac \delta 2, 0)=0,
\label{8}
\end{equation}
i.e.
\begin{equation}
k:=
\frac \alpha {  (T+a) \log(1+\frac Ta) }
=\frac{
 \alpha (8-\delta^2) }
{ 8  \log( 8\delta^{-2}) } \frac 1  T.
\label{9}
\end{equation}
Clearly,
\begin{equation}
k(T+a)=\frac \alpha{  \log( 8\delta^{-2})  }<\alpha.
\label{10}
\end{equation}

Next we choose $D$ so that
\begin{equation}
h(0, 0)=1,
\label{10extra}
\end{equation}
i.e.
$$
D= a^{-k}-F.
$$

Since
$$
h(x, T)=0,\quad \mbox{for}\ |x|=1,
$$
and, in view of (\ref{10}),
$$
D\partial_t h(x, t)=
\frac E { (t+a)^2 }
[\alpha -k(t+a)]> 0,\quad \mbox{for}\ |x|=1,
\ 0\le t\le T,
$$
we have
\begin{equation}
h(x,t)\le 0, \quad \mbox{for}\ |x|=1,
\ 0\le t\le T.
\label{9extra}
\end{equation}

We see from  (\ref{8}), (\ref{10extra})
and the expression of $h$,
that 
$$
h(x, 0)\le h(\frac \delta 2, 0)=0,\quad
\mbox{for}\ |x|\ge \frac \delta 2,
$$
and
$$
h(x, 0))\le h(0,0)
=1,\quad
\mbox{for all}\ x.
$$
Thus, in view of (\ref{A2-parabolic}) and 
(\ref{A3-parabolic}),
\begin{equation}
h\le u\qquad\mbox{on}\
P\partial (B_1\times (0, T]).
\label{11}
\end{equation}

\noindent{\bf Claim.}\ There exists  constants $\overline T>0$,
which depends only on $n, \lambda, \Lambda, \delta$, such that
for all $0<T<\overline T$,
\begin{equation}
(L-\partial_t)h\ge \epsilon,\quad
\mbox{in}\ 
B_1\times (0, T],
\label{12}
\end{equation}
where $\epsilon>0$ is some constant depending only on
$n, \lambda, \Lambda, \delta$, and $T$.

\noindent{\bf Proof of the Claim.}\
We compute
$$
h_i=-\frac {2\alpha x_i}{ t+a} \frac ED,\qquad
h_{ij}=\left( \frac {4\alpha^2 x_ix_j}{  (t+a)^2 }
-\frac {2\alpha \delta_{ij}  }
{ t+a} \right)  \frac ED,
$$
$$
-h_t= \left( \frac k{ t+a} -
\frac{ \alpha|x|^2 }{ (t+a)^2 }\right)  \frac ED.
$$
Thus
\begin{eqnarray*}
J:= 
(t+a) \frac DE (Lh-h_t)
&=& 
\frac { a_{ij} 4\alpha^2 x_ix_j}{ t+a}
-2\alpha \sum_i a_{ii}
-2\alpha b_i x_i+c(t+a)\\
&& -c\frac FE (t+a) +k-\frac{ \alpha |x|^2}{ t+a}.
\end{eqnarray*}
By our choice of $\alpha=1/ \lambda$, and also
$|x|\le 1$, and $c<0$, we have, for some constant
$C$ depending only on $n, \lambda, \Lambda$ and $\delta$,
$$
J\ge k- C(1+T)
= \frac{
  (8-\delta^2) }
{ 8 \lambda  \log( 8\delta^{-2}) } \frac 1  T
-C(1+T).
$$
Clearly, there exists some constant $\overline T>0$ which depends only on
$n, \lambda, \Lambda, \delta$, such that for all $0<T<\overline T$,
we have 
$$
J\ge
 \frac{
  (8-\delta^2) }
{ 9 \lambda  \log( 8\delta^{-2}) } \frac 1  T.
$$
On the other hand
$$
(t+a)\frac DE
\le (T+a) \frac {a^{-k} }  
{E}
\le(T+a)^{k+1} a^{-k}
e^{ \frac {|x|^2}{ \lambda(t+a)} } 
\le (T+a)^{k+1} a^{-k}
e^{ \frac 1{ \lambda a} }.
$$
The claim follows immediately from the above.

\medskip

Let $\overline T>0$  be the  positive constant in the above claim, and
 assume that $T\in (0, \overline T)$.  
Let $\epsilon>0$ be the constant in the claim, which depends on
$T$ in particular, 
and let $u$ satisfy the hypotheses of
Theorem \ref{maximum-parabolic}
with this $\epsilon$.
We will show that
\begin{equation}
u\ge h\qquad\mbox{in}\
 B_1\times (0, T].
\label{13}
\end{equation}
This implies 
$$
u(x, T)\ge \mu(1-|x|), \qquad
\forall\ x\in B_1,
$$
where $\mu>0$ is some constant depending only on
$n, \lambda, \Lambda, \delta$ and $T$.

Since the positive constants $\epsilon$ and $\mu$ 
can clearly be chosen to depend on $T$ monotonically,
the above implies  (\ref{A4-parabolic}).

Now we prove (\ref{13}):
Lower the graph of $h$ to be below that of $u$, in
$B_1\times [0, T]$, and then move it up to a position
$$
h-c_0, 
$$
so that its graph touches that of $u$ from below at some point
$(\bar x, \bar t)$.  It suffices to prove that
$c_0\le 0$.  Suppose not,  $c_0>0$.
Then,
$$
|\bar x|<1, \quad 0<\bar t\le T.
$$
Since $Lu-u_t\le \epsilon$ in viscosity sense, we have,
at $(\bar x, \bar t)$,
$$
L(h-c_0)-h_t\le \epsilon.
$$
It follows, since $c<0$,
$$
(L-\partial_t)h(\bar x, \bar t)\le \epsilon +c c_0<\epsilon.
$$
This contradicts (\ref{12}).  
We have proved (\ref{13}).
Theorem \ref{maximum-parabolic} is established.

\vskip 5pt
\hfill $\Box$
\vskip 5pt

\noindent{\bf Proof of Theorem \ref{thm5.2new}.}\
As usual, we always assume, without loss of generality, that
$c(x,t)\le 0$ on $(0, T]$.

If the assumption (\ref{ABC2}) is replaced by
$$
u\ge 0\quad\mbox{in}\
\Omega\times (0, T],
$$
then the conclusion can be deduced from Theorem \ref{maximum-parabolic}
by using arguments  similar to that used in the proof of
Lemma \ref{lemA-0new}.

Since 
$$
(L-\partial_t)(u+\epsilon t)=0,\quad \mbox{in}\ \Omega\times (0, T],
$$
and
$$
u+\epsilon t\ge u\ge 0,\quad \mbox{on}\ 
P\partial(\Omega\times (0, T]),
$$
we have
$$
u+\epsilon t\ge 0,\quad \Omega\times (0, T].
$$
Thus, as mentioned above, the conclusion of Theorem \ref{thm5.2new}
holds for $u+\epsilon t$.  Namely, for some positive constant
$\bar \mu$ depending only on
 $n, \lambda, \Lambda, \Omega,$
the raduls of $B$, $T_1, T$, 
but independent of $\epsilon$,
such that
\begin{equation}
u+\epsilon t\ge \bar \mu\ dist(x, \partial \Omega),\quad
\mbox{on}\ \Omega\times [T_1/2, T].
\label{B1-1}
\end{equation}

Let $d(x)=dist(x, \partial \Omega)$ denote
the distance of $x$ to $ \partial \Omega$, and we work in
$\Omega\setminus \overline \Omega_\delta$,
where
$$
\Omega_\delta:= \{x\in \Omega\ |\ dist(x)
>\delta\}
$$
for small $\delta$.
The value of $\delta$, depending only on $\Omega$, will be fixed below.

For
$
0<\epsilon\le \frac {\bar \mu \delta}{ 2T},
$
we see from (\ref{B1-1}) that
\begin{equation}
u\ge \bar \mu \delta -\epsilon T\ge \frac{\bar \mu\delta }2,
\quad \mbox{on}\
\partial \Omega_\delta\times [T_1/2, T],
\label{B2-2}
\end{equation}
and
\begin{equation}
u\ge -\epsilon T_1/2+ \bar \mu d(x),\quad
\mbox{on}\ 
(\Omega\setminus \overline \Omega_\delta)\times \{T_1/2\}.
\label{B3-1}
\end{equation}

Fix a function $\rho\in C^\infty([T_1/2, \infty))$ satisfying
$\rho(t)=0, t\ge T_1$;
$-T_1/2\le \rho(t)\le 0, 
T_1/2\le t\le T_1$;
$\rho(T_1/2)=-T_1/2$;
$-2\le \rho'(t)\le 0$, $t\ge 0$.
We use comparison
$$
h(x,t):=
\frac{\bar \mu} 4\left( 
 d(x) +\frac{ d(x)^2}{2\delta}\right)+\epsilon\rho(t).
$$

A computation shows
(see e.g. lemma 7.1 in \cite{CLN1}) that for some small positive numbers
 $\delta$ and $a$, depending only on $\Omega$,
we have
$$
L\left(
 d(x) +\frac{ d(x)^2}{2\delta}\right)
\ge \frac {a\lambda}\delta,\quad \mbox{in}\
(\Omega\setminus \overline \Omega_\delta).
$$
Thus, after further  requiring  that 
$
\epsilon< \frac {a\lambda \bar \mu}{12\delta}$, we have
$$
(L-\partial_t)h
\ge  \frac {a\lambda\bar\mu}{4\delta}-2\epsilon\ge \epsilon,
\qquad\mbox{in}\ 
(\Omega\setminus \overline \Omega_\delta)
\times (0, T].
$$

Now we have 
$$
(L-\partial_t)u \le (L-\partial_t)h,
\qquad\mbox{in}\
(\Omega\setminus \overline \Omega_\delta)
\times [T_1/2, T]
$$
and
$$
u\ge h,
\qquad\mbox{on}\
P\partial(\Omega\setminus \overline \Omega_\delta)
\times [T_1/2, T]),
$$
it follows that
$$
u\ge h,\quad \mbox{in} (\Omega\setminus \overline \Omega_\delta)
\times [T_1/2, T]).
$$
In particular
$$
u\ge h= \frac{\bar \mu}2\left(
 d(x) +\frac{ d(x)^2}{2\delta}\right),
\quad \mbox{in} (\Omega\setminus \overline \Omega_\delta)
\times [T_1,  T]).
$$
 Theorem \ref{thm5.2new} is established.

\vskip 5pt
\hfill $\Box$


\vskip 5pt


\begin{thebibliography}{99}
\bibitem{BD} I. Birindelli and F. Demengel, 
Eigenvalue and Dirichlet problem for fully-nonlinear 
operators in non-smooth domains,
  J. Math. Anal. Appl.  352  (2009),   822-835.
\bibitem{CC} L. Caffarelli and X. Cabre,
Fully nonlinear elliptic equations, American Mathematical
Society Colloquium Publications 43, American Mathematical Society,
Providence, RI, 1995.
\bibitem{CLN1} L. Caffarelli, Y.Y. Li and L. Nirenberg,
Some remarks on singular solutions of
 nonlinear elliptic
equations. I,
	Journal of Fixed Point Theory and Applications
5 (2009), 353-395.
\bibitem{N} L. Nirenberg,
A strong maximum principle for parabolic 
equations, Comm. Pure Appl. Math. 6 (1953),
167-177.
\end{thebibliography}
\end{document}